# DISCUSSION PAPER
## EQUI-ENERGY SAMPLER WITH APPLICATIONS IN STATISTICAL INFERENCE AND STATISTICAL MECHANICS[1,2,3]

By S. C. Kou, Qing Zhou and Wing Hung Wong

*Harvard University, Harvard University and Stanford University*

We introduce a new sampling algorithm, the equi-energy sampler, for efficient statistical sampling and estimation. Complementary to the widely used temperature-domain methods, the equi-energy sampler, utilizing the temperature–energy duality, targets the energy directly. The focus on the energy function not only facilitates efficient sampling, but also provides a powerful means for statistical estimation, for example, the calculation of the density of states and microcanonical averages in statistical mechanics. The equi-energy sampler is applied to a variety of problems, including exponential regression in statistics, motif sampling in computational biology and protein folding in biophysics.

**1. Introduction.** Since the arrival of modern computers during World War II, the Monte Carlo method has greatly expanded the scientific horizon to study complicated systems ranging from the early development in computational physics to modern biology. At the heart of the Monte Carlo method lies the difficult problem of sampling and estimation: Given a target distribution, usually multidimensional and multimodal, how do we draw samples from it and estimate the statistical quantities of interest? In this article, we attempt to introduce a new sampling algorithm, the equi-energy sampler, to address the problem. Since the Monte Carlo method began from calculations in statistical physics and mechanics, to introduce the equi-energy sampler, we begin from statistical mechanics.

Received June 2004; revised March 2005.

[1]Supported in part by NSF, NIH and Harvard University Clarke–Cooke Fund.

[2]Discussed in 10.1214/009053606000000470, 10.1214/009053606000000489, 10.1214/009053606000000498 and 10.1214/009053606000000506; rejoinder at 10.1214/009053606000000524.

[3]S. C. Kou and Qing Zhou contributed equally to this work.

The starting point of a statistical mechanical computation is the energy function or Hamiltonian $h(x)$. According to Boltzmann and Gibbs, the distribution of a system in thermal equilibrium at temperature $T$ is described by the Boltzmann distribution,

$$(1) \qquad p_T(x) = \frac{1}{Z(T)} \exp(-h(x)/T),$$

where $Z(T) = \sum_x \exp(-h(x)/T)$ is referred to as the partition function. For any state function $g(x)$, its expectation $\mu_g(T)$ with respect to the Boltzmann distribution is called its Boltzmann average, also known as the thermal average in the physics literature,

$$(2) \qquad \mu_g(T) = \sum_x g(x) \exp(-h(x)/T)/Z(T).$$

To study the system, in many cases we are interested in using Monte Carlo simulation to obtain estimates of Boltzmann averages as functions of temperature for various state functions. In addition to Boltzmann averages, estimating the partition function $Z(T)$, which represents the dependency of the normalization factor in (1) as a function of temperature, is also of significant interest, as it is well known that many important thermodynamic quantities, such as free energy, specific heat, internal energy and so on, can be computed directly from the partition function (see Section 4). The fundamental algorithm for computing Boltzmann averages is due to Metropolis et al. [29], who proposed the use of a reversible Markov chain constructed in such a way so as to guarantee that its stationary distribution is the Boltzmann distribution (1). Later, this algorithm was generalized by Hastings [11] to allow the use of an asymmetric transition kernel. Given a current state $x$, this Metropolis–Hastings algorithm generates a new state by either reusing the current state $x$ or moving to a new state $y$ drawn from a proposal kernel $K(x \to y)$. The proposal state $y$ is accepted with probability $\min(1, MR)$ where $MR$ is the Metropolis–Hastings ratio $p_T(y)K(y \to x)/p_T(x)K(x \to y)$. The algorithm in this way generates a Markov chain $X_i$, $i = 1, \ldots, n$. Under ergodic conditions [39], the time average $n^{-1} \sum_{i=1}^n g(X_i)$ provides a consistent estimate of the Boltzmann average (2).

The Metropolis algorithm, however, can perform poorly if the energy function has many local minima separated by high barriers that cannot be crossed by the proposal moves. In this situation the chain will be trapped in local energy wells and will fail to sample the Boltzmann distribution correctly. To overcome this problem one can design specific moves that have a higher chance to cut across the energy barrier (e.g., the conditional sampling moves in Gibbs sampling) or to add auxiliary variables so that the energy wells become connected by the added dimension (e.g., the group Ising updating of Swendsen and Wang [37], or the data-augmentation technique of Tanner and Wong [38]). However, these remedies are problem-specific and may or may not work for any given problem. A breakthrough occurred with the development of the parallel tempering algorithm by Geyer [8] (also called exchanged Monte Carlo; Hukushima and Nemoto [13]). The idea is to perform parallel Metropolis sampling at different temperatures. Occasionally one proposes to exchange the states of two neighboring chains (i.e., chains with adjacent temperature levels). The acceptance probability for the exchange is designed to ensure that the joint states of all the parallel chains evolve according to the Metropolis–Hastings rule with the product distribution (i.e., the product of the Boltzmann distributions at the different temperatures) as the target distribution. Geyer's initial objective was to use the (hopefully) faster mixing of the high-temperature chains to drive the mixing of the whole system, thereby to achieve faster mixing at the low-temperature chain as well. It is clear that parallel tempering can provide estimates of Boltzmann averages at the temperatures used in the simulation. Marinari and Parisi [27] developed simulated tempering that uses just a single chain but augments the state by a temperature variable that is dynamically moved up or down the temperature ladder. These authors further developed the theory of using the samples in the multiple temperatures to construct estimates of the partition function, and to investigate phase transition. In the meantime, Geyer [9] also proposed a maximum likelihood approach to estimate the ratios of normalization constants and hence obtain information on the partition function at the selected temperatures.

In contrast to the statistical mechanical computations, the starting point in statistical inference is usually *one* given distribution, for example, the distribution on a high-dimensional parameter space. If we take the energy to be the negative log-density function in this case, we are then interested in obtaining the Boltzmann average only at $T = 1$. The Metropolis and related algorithms have been developed and applied to solve many statistical computation problems, and have greatly enhanced our ability to analyze problems ranging from image analysis to missing value problems to biological sequence analysis to single-molecule chemistry [6, 7, 17, 20, 22, 38]. However, as Geyer, Marinari, Parisi and others have pointed out, even if the immediate interest is at $T = 1$, simulation at temperatures other than $T = 1$ is often necessary in order to achieve efficient sampling. Furthermore, computing the normalization constants (i.e., partition function) is also important in statistical tasks such as the determination of likelihood ratios and Bayes factors [9, 10, 28].

It is thus seen that historically dynamic Monte Carlo methods were developed to simulate from the Boltzmann distribution at fixed temperatures. These methods aim to provide direct estimates of parameters such as Boltzmann averages and partition functions, which are functions of temperature.

We hence refer to these methods as *temperature-domain* methods. The purpose of this article is to develop an *alternative* sampling and estimation approach based on *energy-domain* considerations. We will construct algorithms for the direct estimation of parameters such as microcanonical averages and density of states (see Section 2) that are functions of *energy*. We will see in Section 2 that there is a duality between temperature-domain functions and energy-domain functions, so that once we have obtained estimates of the density of states and microcanonical averages (both are energy-domain functions), we can easily transfer to the temperature domain to obtain the partition function and the Boltzmann averages. In Section 3 we introduce the equi-energy sampler (EE sampler), which, targeting energy directly, is a new Monte Carlo algorithm for the efficient sampling from multiple energy intervals. In Sections 4 and 5 we explain how to use these samples to obtain estimates of density of states and microcanonical averages, and how to extend the energy-domain method to estimate statistical quantities in general. In Section 6 we illustrate the wide applicability of this method by applying the equi-energy sampler and the estimation methods to a variety of problems, including an exponential regression problem, the analysis of regulatory DNA motifs and the study of a simplified model for protein folding. Section 7 concludes the article with discussion and further remarks.

**2. Energy–temperature duality.** The Boltzmann law (1) implies that the conditional distribution of the system given its energy $h(x) = u$ is the uniform distribution on the equi-energy surface $\{x : h(x) = u\}$. In statistical mechanics, this conditional distribution is referred to as the microcanonical distribution given energy $u$. Accordingly, the conditional expectation of a state function $g(x)$ given an energy level $u$ is called its microcanonical average:

$$(3) \qquad \nu_g(u) = E(g(X)|h(X) = u).$$

Note that (3) is independent of the temperature $T$ used in the Boltzmann distribution for $X$. Suppose that the infinitesimal volume of the energy slice $\{x : h(x) \in (u, u + du)\}$ is approximately equal to $\Omega(u) \, du$. This function $\Omega(u)$ is then called the *density of states* function. If the state space is discrete, then we replace the volume by counts, in which case $\Omega(u)$ is simply the number of states with the energy equal to $u$. Without loss of generality, we assume that the minimum energy of the system $u_{\min} = 0$. The following result follows easily from these definitions.

LEMMA 1. *Let $\beta = 1/T$ denote the inverse temperature so that the Boltzmann averages and partition function are indexed by $\beta$ as well as by $T$; then*

$$\mu_g(\beta^{-1})Z(\beta^{-1}) = \int_0^\infty \nu_g(u)\Omega(u)e^{-\beta u} \, du.$$

*In particular, the partition function $Z(\beta^{-1})$ and the density of states $\Omega(u)$ form a Laplace transform pair.*

This lemma suggests that the Boltzmann averages and the partition function can be obtained through Monte Carlo algorithms designed to compute the density of states and microcanonical averages. We hence refer to such algorithms as energy-domain algorithms.

The earliest energy-domain Monte Carlo algorithm is the multicanonical algorithm due to Berg and Neuhaus [2], which aims to sample from a distribution flat in the energy domain through an iterative estimation updating scheme. Later, the idea of iteratively updating the target distribution was generalized to histogram methods (see [18, 40] for a review). The main purpose of these algorithms is to obtain the density of states and related functions such as the specific heat. They do not directly address the estimation of Boltzmann averages.

In this article we present a different method that combines the use of multiple energy ranges, multiple temperatures and step sizes, to produce an efficient sampling scheme capable of providing direct estimates of all microcanonical averages as well as the density of states. We do not use iterative estimation of density of states as in the multicanonical approach; instead, the key of our algorithm is a new type of move called the equi-energy jump that aims to move directly between states with similar energy (see the next section). The relationship between the multicanonical algorithm and the equi-energy sampler will be discussed further in Section 7.

## 3. The equi-energy sampler.

3.1. *The algorithm.* In Monte Carlo statistical inference one crucial task is to obtain samples from a given distribution, often known up to a normalizing constant. Let $\pi(x)$ denote the target distribution and let $h(x)$ be the associated energy function. Then $\pi(x) \propto \exp(-h(x))$. For simple problems, the famous Metropolis–Hastings (MH) algorithm, which employs a local Markov chain move, could work. However, if $\pi(x)$ is multimodal and the modes are far away from each other, which is often the case for practical multidimensional distributions, algorithms relying on local moves such as the MH algorithm or the Gibbs sampler can be easily trapped in a local mode indefinitely, resulting in inefficient and even unreliable samples.

The EE sampler aims to overcome this difficulty by working on the energy function directly. First, a sequence of energy levels is introduced:

(4) $$H_0 < H_1 < H_2 < \cdots < H_K < H_{K+1} = \infty,$$

such that $H_0$ is below the minimum energy, $H_0 \leq \inf_x h(x)$. Associated with the energy levels is a sequence of temperatures

$$1 = T_0 < T_1 < \cdots < T_K.$$

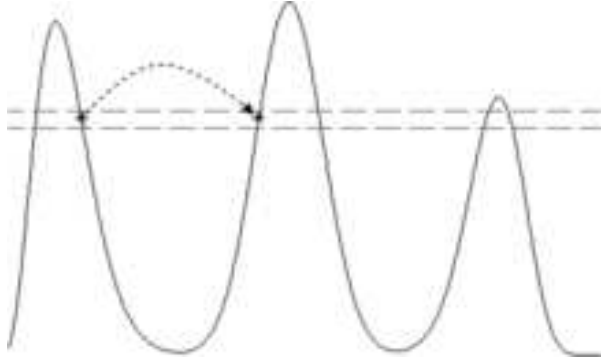

FIG. 1.  *Illustration of the equi-energy jump, where the sampler can jump freely between the states with similar energy levels.*

The EE sampler considers $K + 1$ distributions, each indexed by a temperature and an energy truncation. The energy function of the $i$th distribution $\pi_i$ $(0 \leq i \leq K)$ is $h_i(x) = \frac{1}{T_i}(h(x) \vee H_i)$, that is, $\pi_i(x) \propto \exp(-h_i(x))$. For each $i$, a sampling chain targeting $\pi_i$ is constructed. Clearly $\pi_0$ is the initial distribution of interest. The EE sampler employs the other $K$ chains to overcome local trapping, because for large $i$ the energy truncation and the high temperature on $h_i(x)$ flatten the distribution $\pi_i(x)$, making it easier to move between local modes. The quick mixing of chains with large $i$ is utilized by the EE sampler, through a step termed the *equi-energy jump*, to help sampling from $\pi_i$ with small $i$, where the landscape is more rugged.

The equi-energy jump, illustrated in Figure 1, aims to directly move between states with similar energy levels. Intuitively, if it can be implemented properly in a sampling algorithm, it will effectively eliminate the problem of local trap. The fact that the $i$th energy function $h_i(x)$ is monotone in $h(x)$ implies that, above the truncation values, the equi-energy sets $S_i(H) = \{x : h_i(x) = H\}$ are mutually consistent across $i$. Thus once we have constructed an empirical version of the equi-energy sets at high-order $\pi_i$ (i.e., $\pi_i$ with large $i$), these high-order empirical sets will remain valid at low-order $\pi_i$. Therefore, after the empirical equi-energy sets are first constructed at high order $\pi_i$, the local trapping at low-order $\pi_i$ can be largely evaded by performing an equi-energy jump that allows the current state to jump to another state drawn from the already constructed high-order empirical equi-energy set that has energy level close to the current state. This is the basic idea behind the EE sampler. We will refer to the empirical equi-energy sets as energy rings hereafter.

To construct the energy rings, the state space $\mathcal{X}$ is partitioned according to the energy levels, $\mathcal{X} = \bigcup_{j=0}^{K} D_j$, where $D_j = \{x : h(x) \in [H_j, H_{j+1})\}$, $0 \leq j \leq K$, are the energy sets, determined by the energy sequence (4). For any $x \in \mathcal{X}$, let $I(x)$ denote the partition index such that $I(x) = j$, if $x \in D_j$, that is, if $h(x) \in [H_j, H_{j+1})$.

The EE sampler begins from an MH chain $X^{(K)}$ targeting the highest-order distribution $\pi_K$. After an initial burn-in period, the EE sampler starts constructing the $K$th-order energy rings $\hat{D}_j^{(K)}$ by grouping the samples according to their energy levels; that is, $\hat{D}_j^{(K)}$ consists of all the samples $X_n^{(K)}$ such that $I(X_n^{(K)}) = j$. After the chain $X^{(K)}$ has been running for $N$ steps, the EE sampler starts the second highest-order chain $X^{(K-1)}$ targeting $\pi_{K-1}$, while it keeps on running $X^{(K)}$ and updating $\hat{D}_j^{(K)}$ ($0 \le j \le K$). The chain $X^{(K-1)}$ is updated through two operations: the local move and the equi-energy jump. At each update a coin is flipped; with probability $1 - p_{\text{ee}}$ the current state $X_n^{(K-1)}$ undergoes an MH local move to give the next state $X_{n+1}^{(K-1)}$, and with probability $p_{\text{ee}}$, $X_n^{(K-1)}$ goes through an equi-energy jump. In the equi-energy jump, a state $y$ is chosen uniformly from the highest-order energy ring $\hat{D}_j^{(K)}$ indexed by $j = I(X_n^{(K-1)})$ that corresponds to the energy level of $X_n^{(K-1)}$ [note that $y$ and $X_n^{(K-1)}$ have similar energy level, since $I(y) = I(X_n^{(K-1)})$ by construction]; the chosen $y$ is accepted to be the next state $X_{n+1}^{(K-1)}$ with probability $\min(1, \frac{\pi_{K-1}(y)\pi_K(X_n^{(K-1)})}{\pi_{K-1}(X_n^{(K-1)})\pi_K(y)})$; if $y$ is not accepted, $X_{n+1}^{(K-1)}$ keeps the old value $X_n^{(K-1)}$. After a burn-in period on $X^{(K-1)}$, the EE sampler starts the construction of the second highest-order [i.e., $(K-1)$st order] energy rings $\hat{D}_j^{(K-1)}$ in much the same way as the construction of $\hat{D}_j^{(K)}$, that is, collecting the samples according to their energy levels. Once the chain $X^{(K-1)}$ has been running for $N$ steps, the EE sampler starts $X^{(K-2)}$ targeting $\pi_{K-2}$ while it keeps on running $X^{(K-1)}$ and $X^{(K)}$. Like $X^{(K-1)}$, the chain $X^{(K-2)}$ is updated by the local MH move and the equi-energy jump with probabilities $1 - p_{\text{ee}}$ and $p_{\text{ee}}$, respectively. In the equi-energy jump, a state $y$ uniformly chosen from $\hat{D}_{I(X_n^{(K-2)})}^{(K-1)}$, where $X_n^{(K-2)}$ is the current state, is accepted to be the next state $X_{n+1}^{(K-2)}$ with probability $\min(1, \frac{\pi_{K-2}(y)\pi_{K-1}(X_n^{(K-2)})}{\pi_{K-2}(X_n^{(K-2)})\pi_{K-1}(y)})$. The EE sampler thus successively moves down the energy and temperature ladder until the distribution $\pi_0$ is reached. Each chain $X^{(i)}$, $0 \le i < K$, is updated by the equi-energy jump and the local MH move; the equi-energy move proposes a state $y$ uniformly chosen from the energy ring $\hat{D}_{I(X_n^{(i)})}^{(i+1)}$ and accepts the proposal with probability $\min(1, \frac{\pi_i(y)\pi_{i+1}(X_n^{(i)})}{\pi_i(X_n^{(i)})\pi_{i+1}(y)})$. At each chain $X^{(i)}$, the energy rings $\hat{D}_j^{(i)}$ are constructed after a burn-in period, and will be used for chain $X^{(i-1)}$ in the equi-energy jump. The basic sampling scheme can be summarized as follows.

**The algorithm of the EE sampler**

Assign $X_0^{(i)}$ an initial value and set $\hat{D}_j^{(i)} = \varnothing$ for all $i$ and $j$

For $n = 1, 2, \ldots$

    For $i = K$ downto 0

      if $n > (K - i)(B + N)$ do

      (* $B$ is the burn-in period, and $N$ is the

          period to construct initial energy rings *)

        if $i = K$ or if $\hat{D}_{I(X_{n-1}^{(i)})}^{(i+1)} = \varnothing$ do

          perform an MH step on $X_{n-1}^{(i)}$ with

             target distribution $\pi_i$ to obtain $X_n^{(i)}$

        else

          with probability $1 - p_{\mathrm{ee}}$, perform an MH step

             on $X_{n-1}^{(i)}$ targeting $\pi_i$ to obtain $X_n^{(i)}$

          with probability $p_{\mathrm{ee}}$, uniformly pick a state $y$ from $\hat{D}_{I(X_{n-1}^{(i)})}^{(i+1)}$ and

let

$$X_n^{(i)} \leftarrow y \text{ with prob } \min\left(1, \frac{\pi_i(y)\pi_{i+1}(X_{n-1}^{(i)})}{\pi_i(X_{n-1}^{(i)})\pi_{i+1}(y)}\right);$$

$$X_n^{(i)} \leftarrow X_{n-1}^{(i)} \text{ with the remaining prob.}$$

        endif

        if $n > (K - i)(B + N) + B$ do

          $\hat{D}_{I(X_n^{(i)})}^{(i)} \leftarrow \hat{D}_{I(X_n^{(i)})}^{(i)} + \{X_n^{(i)}\}$

          (* add the sample to the energy rings after the burn-in period

*)

        endif

      endif

    endfor

  endfor

The idea of moving along the equi-energy surface bears a resemblance to the auxiliary variable approach, where to sample from a target distribution $\pi(x)$, one can iteratively first sample an auxiliary variable $U \sim \mathrm{Unif}[0, \pi(X)]$, and then sample $X \sim \mathrm{Unif}\{x : \pi(x) \geq U\}$. This approach has been used by Edwards and Sokal [5] to explain the Swendsen–Wang [37] clustering algorithm for Ising model simulation, and by Besag and Green [3] and Higdon [12] on spatial Bayesian computation and image analysis. Roberts and Rosenthal [32], Mira, Moller and Roberts [30] and Neal [31] provide further discussion under the name of slice sampling. Under our setting, this auxiliary variable (slice sampler) approach amounts to sampling from the lower-energy sets. In comparison, the EE sampler's focus on the equi-energy set is motivated directly from the concept of microcanonical distribution in statistical mechanics. More importantly, the equi-energy jump step in the EE sampler offers a *practical* means to carry out the idea of moving along the energy sets (or moving horizontally along the density contours), and thus provides an effective way to study not only systems in statistical mechanics, but also general statistical inference problems (see Sections 4, 5 and 6).

3.2. *The steady-state distribution.* With help from the equi-energy jump to address the local trapping, the EE sampler aims to efficiently draw samples from the given distribution $\pi$ (which is identical to $\pi_0$). A natural question is then: In the long run, will the EE samples follow the correct distribution?

The following theorem shows that the steady-state distribution of chain $X^{(i)}$ is indeed $\pi_i$; in particular, the steady-state distribution of $X^{(0)}$ is $\pi_0 = \pi$.

THEOREM 2. *Suppose* (i) *the highest-order chain* $X^{(K)}$ *is irreducible and aperiodic,* (ii) *for* $i = 0, 1, \ldots, K - 1$, *the MH transition kernel* $T_{\mathrm{MH}}^{(i)}$ *of* $X^{(i)}$ *connects adjacent energy sets in the sense that for any* $j$ *there exist sets* $A_1 \subset D_j$, $A_2 \subset D_j$, $B_1 \subset D_{j-1}$ *and* $B_2 \subset D_{j+1}$ *with positive measure such that the transition probabilities*

$$T_{\mathrm{MH}}^{(i)}(A_1, B_1) > 0, \qquad T_{\mathrm{MH}}^{(i)}(A_2, B_2) > 0$$

*and* (iii) *the energy set probabilities* $p_j^{(i)} = P_{\pi_i}(X \in D_j) > 0$ *for all* $i$ *and* $j$. *Then* $X^{(i)}$ *is ergodic with* $\pi_i$ *as its steady-state distribution.*

PROOF. We use backward induction to prove the theorem.

For $i = K$, $X^{(K)}$ simply follows the standard MH scheme. The desired conclusion thus follows from the fact that $X^{(K)}$ is aperiodic and irreducible.

Now assume the conclusion holds for the $(i+1)$st order chain, that is, assume $X^{(i+1)}$ is ergodic with steady-state distribution $\pi_{i+1}$. We want to show that the conclusion also holds for $X^{(i)}$. According to the construction of the EE sampler, if at the $n$th step $X_n^{(i)} = x$, then at the next step with probability $1 - p_{\mathrm{ee}}$ $X_{n+1}^{(i)}$ will be drawn from the transition kernel $T_{\mathrm{MH}}^{(i)}(x, \cdot)$, and with probability $p_{\mathrm{ee}}$ $X_{n+1}^{(i)}$ will be equal to a $y$ from $\hat{D}_{I(x)}^{(i+1)}$ with probability

$$P(X_{n+1}^{(i)} = y) = \frac{1}{|\hat{D}_{I(x)}^{(i+1)}|} \min\left(1, \frac{\pi_i(y)\pi_{i+1}(x)}{\pi_i(x)\pi_{i+1}(y)}\right), \qquad y \in \hat{D}_{I(x)}^{(i+1)}.$$

Therefore for any measurable set $A$, the conditional probability

$$P(X_{n+1}^{(i)} \in A | X_n^{(i)} = x, X^{(i+1)})$$

$$= (1 - p_{\text{ee}}) T_{\text{MH}}^{(i)}(x, A)$$

$$+ p_{\text{ee}} \frac{1}{|\hat{D}_{I(x)}^{(i+1)}|} \sum_{y \in \hat{D}_{I(x)}^{(i+1)}} I(y \in A) \min\left(1, \frac{\pi_i(y)\pi_{i+1}(x)}{\pi_i(x)\pi_{i+1}(y)}\right)$$

$$+ p_{\text{ee}} \left[1 - \frac{1}{|\hat{D}_{I(x)}^{(i+1)}|} \sum_{y \in \hat{D}_{I(x)}^{(i+1)}} \min\left(1, \frac{\pi_i(y)\pi_{i+1}(x)}{\pi_i(x)\pi_{i+1}(y)}\right)\right] I(x \in A).$$

Using the induction assumption of the ergodicity of $X^{(i+1)}$ and also the fact that the lower-order chain $X^{(i)}$ does not affect the higher-order chain $X^{(i+1)}$, we have, as $n \to \infty$,

$$P(X_{n+1}^{(i)} \in A | X_n^{(i)} = x)$$

$$= \int P(X_{n+1}^{(i)} \in A | X_n^{(i)} = x, X^{(i+1)}) \, dP(X^{(i+1)} | X_n^{(i)} = x)$$

$$= \int P(X_{n+1}^{(i)} \in A | X_n^{(i)} = x, X^{(i+1)}) \, dP(X^{(i+1)})$$

(5)

$$\to (1 - p_{\text{ee}}) T_{\text{MH}}^{(i)}(x, A)$$

$$+ p_{\text{ee}} \frac{1}{p_{I(x)}^{(i+1)}} \int_{y \in A \cap D_{I(x)}} \pi_{i+1}(y) \min\left(1, \frac{\pi_i(y)\pi_{i+1}(x)}{\pi_i(x)\pi_{i+1}(y)}\right) dy$$

$$+ p_{\text{ee}} \left[1 - \frac{1}{p_{I(x)}^{(i+1)}} \int_{y \in D_{I(x)}} \pi_{i+1}(y) \min\left(1, \frac{\pi_i(y)\pi_{i+1}(x)}{\pi_i(x)\pi_{i+1}(y)}\right) dy\right] I(x \in A).$$

Similarly, as $n \to \infty$, the difference

$$P(X_{n+1}^{(i)} \in A | X_n^{(i)} = x, X_{n-1}^{(i)}, \ldots, X_1^{(i)}) - P(X_{n+1}^{(i)} \in A | X_n^{(i)} = x) \to 0.$$

Now let us define a new transition kernel $S^{(i)}(x, \cdot)$, which undergoes the transition $T_{\text{MH}}^{(i)}(x, \cdot)$ with probability $1 - p_{\text{ee}}$, and with probability $p_{\text{ee}}$ undergoes an MH transition with the proposal density $q(x, y) = \frac{1}{p_{I(x)}^{(i+1)}} \pi_{i+1}(y) I(y \in D_{I(x)})$, that is, $\pi_{i+1}(y)$ confined to the energy set $D_{I(x)}$. We then note that the right-hand side of (5) corresponds exactly to the transition kernel $S^{(i)}(x, \cdot)$. Therefore, under the induction assumption, $X^{(i)}$ is asymptotically equivalent to a Markovian sequence governed by $S^{(i)}(x, \cdot)$.

Since the kernel $T_{\text{MH}}^{(i)}(x, \cdot)$ connects adjacent energy sets and the proposal $q(x, y)$ connects points in the same equi-energy set, it follows from Chapman–Kolmogorov and $0 < p_{\text{ee}} < 1$ that $S^{(i)}(x, \cdot)$ is irreducible. $S^{(i)}(x, \cdot)$ is also aperiodic because the proposal $q(x, y)$ has positive probability to leave the configuration $x$ staying the same.

Since $S^{(i)}(x, \cdot)$ keeps $\pi_i$ as the steady-state distribution, it finally follows from the standard Markov chain convergence theorem and the asymptotic equivalence (5) that $X^{(i)}$ is ergodic with $\pi_i$ as its steady-state distribution. The proof is thus terminated. □

REMARK 3. The assumption (ii) is weaker than assuming that $T_{\text{MH}}^{(i)}$ is irreducible for $i = 0, 1, \ldots, K - 1$, because we can see that essentially the function of the MH local move is to bridge adjacent energy sets, while the equi-energy jump allows jumps within an equi-energy set.

3.3. *Practical implementation.* There are some flexibilities in the practical implementation of the EE sampler. We provide some suggestions based on our own experience.

1. The choice of the temperature and energy ladder.

   Given the lowest and second highest energy levels $H_0$ and $H_K$, we found that setting the other energy levels by a geometric progression, or equivalently setting $\log(H_{i+1} - H_i)$ to be evenly spaced, often works quite well. The temperature could be chosen such that $(H_{i+1} - H_i)/T_i \approx c$, and we found that $c \in [1, 5]$ often works well.

2. The choice of $K$, the number of temperature and energy levels.

   The choice of $K$ depends on the complexity of the problem. More chains and energy levels are usually needed if the target distribution is high-dimensional and multimodal. In our experience $K$ could be roughly proportional to the dimensionality of the target distribution.

3. The equi-energy jump probability $p_{\text{ee}}$.

   In our experience taking $p_{\text{ee}} \in [5\%, 30\%]$ often works quite well. See Section 3.4 for more discussion.

4. Self-adaptation of the MH-proposal step size.

   As the order $i$ increases, the distribution $\pi_i$ becomes more and more flat. Intuitively, to efficiently explore a flat distribution, one should use a large step size in the MH proposal, whereas for a rough distribution, the step size has to be small. Therefore, in the EE sampler each chain $X^{(i)}$ should have its own step size in the local MH exploration. In practice, however, it is often difficult to choose the right step sizes in the very beginning. One can hence let the sampler tune by itself the step sizes. For each chain, the sampler can from time to time monitor the acceptance rate in the MH local move, and increase (decrease) the step size by a fixed factor, if the acceptance rate is too high (low). Note that in this self-adaptation the energy-ring structure remains unchanged.

5. **Adjusting the energy and temperature ladder.**

In many problems, finding a close lower bound $H_0$ for the energy function $h(x)$ is not very difficult. But in some cases, especially when $h(x)$ is difficult to optimize, one might find during the sampling that the energy value at some state is actually smaller than the pre-assumed lower bound $H_0$. If this happens, we need to adjust the energy ladder and the temperatures, because otherwise the energy sets $D_j$ would not have the proper sizes to cover the state space, which could affect the sampling efficiency. The adjustment can be done by dynamically monitoring the sampler. Suppose after the $i$th chain $X^{(i)}$ is started, but before the $(i-1)$st chain gets started, we find that the lowest energy value $H_{\min}$ reached so far is smaller than $H_0$. Then we first reset $H_0 = H_{\min} - \beta$, where the constant $\beta > 0$, say $\beta = 2$. Next given $H_i$ and the new $H_0$ we reset the in-between energy levels by a geometric progression, and if necessary add in more energy levels between $H_0$ and $H_i$ (thus adding more chains) so that the sequence $H_{j+1} - H_j$ is still monotone increasing in $j$. The temperatures between $T_0 = 1$ and $T_i$ are reset by $(H_{j+1} - H_j)/T_j \approx c$. With the energy ladder adjusted, the samples are regrouped to new energy rings. Note that since the chains $X^{(K)}, X^{(K-1)}, \ldots, X^{(i)}$ have already started, we do not change the values of $H_K, \ldots, H_i$ and $T_K, \ldots, T_i$, so that the target distributions $\pi_K, \ldots, \pi_i$ are not altered.

3.4. *A multimodal illustration.* As an illustration, we consider sampling from a two-dimensional normal mixture model taken from [23],

$$(6) \qquad f(\mathbf{x}) = \sum_{i=1}^{20} \frac{w_i}{2\pi\sigma_i^2} \exp\left\{-\frac{1}{2\sigma_i^2}(\mathbf{x} - \boldsymbol{\mu}_i)'(\mathbf{x} - \boldsymbol{\mu}_i)\right\},$$

where $\sigma_1 = \cdots = \sigma_{20} = 0.1$, $w_1 = \cdots = w_{20} = 0.05$, and the 20 mean vectors

$$(\boldsymbol{\mu}_1, \boldsymbol{\mu}_2, \ldots, \boldsymbol{\mu}_{20}) = \left(\begin{matrix} 2.18 & 8.67 & 4.24 & 8.41 & 3.93 & 3.25 & 1.70 \\ 5.76 & 9.59 & 8.48 & 1.68 & 8.82 & 3.47 & 0.50 \end{matrix}\right.,$$

$$\begin{matrix} 4.59 & 6.91 & 6.87 & 5.41 & 2.70 & 4.98 & 1.14 \\ 5.60 & 5.81 & 5.40 & 2.65 & 7.88 & 3.70 & 2.39 \end{matrix},$$

$$\left.\begin{matrix} 8.33 & 4.93 & 1.83 & 2.26 & 5.54 & 1.69 \\ 9.50 & 1.50 & 0.09 & 0.31 & 6.86 & 8.11 \end{matrix}\right).$$

Since most local modes are more than 15 standard deviations away from the nearest ones [see Figure 2(a)], this mixture distribution poses a serious challenge for sampling algorithms, and thus serves as a good test. We applied the EE sampler to this problem. Since the minimum value of the energy function $h(\mathbf{x}) = -\log(f(\mathbf{x}))$ is around $-\log(\frac{5}{2\pi}) = 0.228$, we took $H_0 = 0.2$. $K$ was set to 4, so only five chains were employed. The energy levels $H_1, \ldots, H_4$ were set by a geometric progression in the interval $[0, 200]$. The settings for energy levels and temperature ladders are tabulated in Table 1. The equi-energy jump probability $p_{\text{ee}}$ was taken to be 0.1. The initial states of the chains $X^{(i)}$ were drawn uniformly from $[0, 1]^2$, a region far from the centers $\boldsymbol{\mu}_1, \boldsymbol{\mu}_2, \ldots, \boldsymbol{\mu}_{20}$, so as to make the sampling challenging. The MH proposal is taken to be bivariate Gaussian: $X_{n+1}^{(i)} \sim N_2(X_n^{(i)}, \tau_i^2 I_2)$, where the initial MH proposal step size $\tau_i$ for the $i$th order chain $X^{(i)}$ was taken to be $0.25\sqrt{T_i}$. The step size was finely tuned later in the algorithm such that the acceptance ratio was in the range $(0.22, 0.32)$. After a burn-in period, each chain was run for 50,000 iterations. Figure 2 shows the samples generated in each chain: With the help of the higher-order chains, where the distributions are more flat, all the modes of the target distribution were successfully visited by $X^{(0)}$. The number of samples in each energy ring is reported in Table 1. One can see that for low-order chains the samples are mostly concentrated in the low-energy rings, while for high-order chains more samples are distributed in the high-energy rings.

As a comparison, we also applied parallel tempering (PT) [8] to this problem. The PT procedure also adopts a temperature ladder; it uses a swap between neighboring temperature chains to help the low-temperature chain move. We ran the PT to sample from (6) with the same parameter and initialization setting. The step size of PT was tuned to make the acceptance ratio of the MH move between 0.22 and 0.32. The exchange (swap) probability of PT was taken to be 0.1 to make it comparable with $p_{\text{ee}} = 0.1$ in the EE sampler, and in each PT exchange operation, $K = 4$ swaps were proposed to exchange samples in neighboring chains. The overall acceptance rates for the MH move in the EE sampler and parallel tempering were 0.27 and 0.29, respectively. In the EE sampler, the acceptance rate for the equi-energy jump was 0.82, while the acceptance rate for the exchange operation in PT was 0.59. Figure 3(a) shows the path of the last 2000 samples in $X^{(0)}$ for the EE sampler, which visited all the 20 components frequently; in comparison PT only visited 14 components in the same number of samples [Figure 3(b)]. As a further test, we ran the EE sampler and the PT 20 times independently, and sought to estimate the mean vector $(EX_1, EX_2)$ and the second moment $(EX_1^2, EX_2^2)$ using the samples generated from the target chain $X^{(0)}$. The results are shown in Table 2 (the upper half). It is clear that the EE sampler provided more accurate estimates with smaller mean squared errors.

TABLE 1
*Sample size of each energy ring*

| Chain | Energy rings | | | | |
|---|---|---|---|---|---|
| | < 2.0 | [2.0, 6.3) | [6.3, 20.0) | [20.0, 63.2) | ≥ 63.2 |
| $X^{(0)}, T_0 = 1$ | 41631 | 8229 | 140 | 0 | 0 |
| $X^{(1)}, T_1 = 2.8$ | 21118 | 23035 | 5797 | 50 | 0 |
| $X^{(2)}, T_2 = 7.7$ | 7686 | 16285 | 22095 | 3914 | 20 |
| $X^{(3)}, T_3 = 21.6$ | 3055 | 6470 | 17841 | 20597 | 2037 |
| $X^{(4)}, T_4 = 60.0$ | 1300 | 2956 | 8638 | 20992 | 16114 |

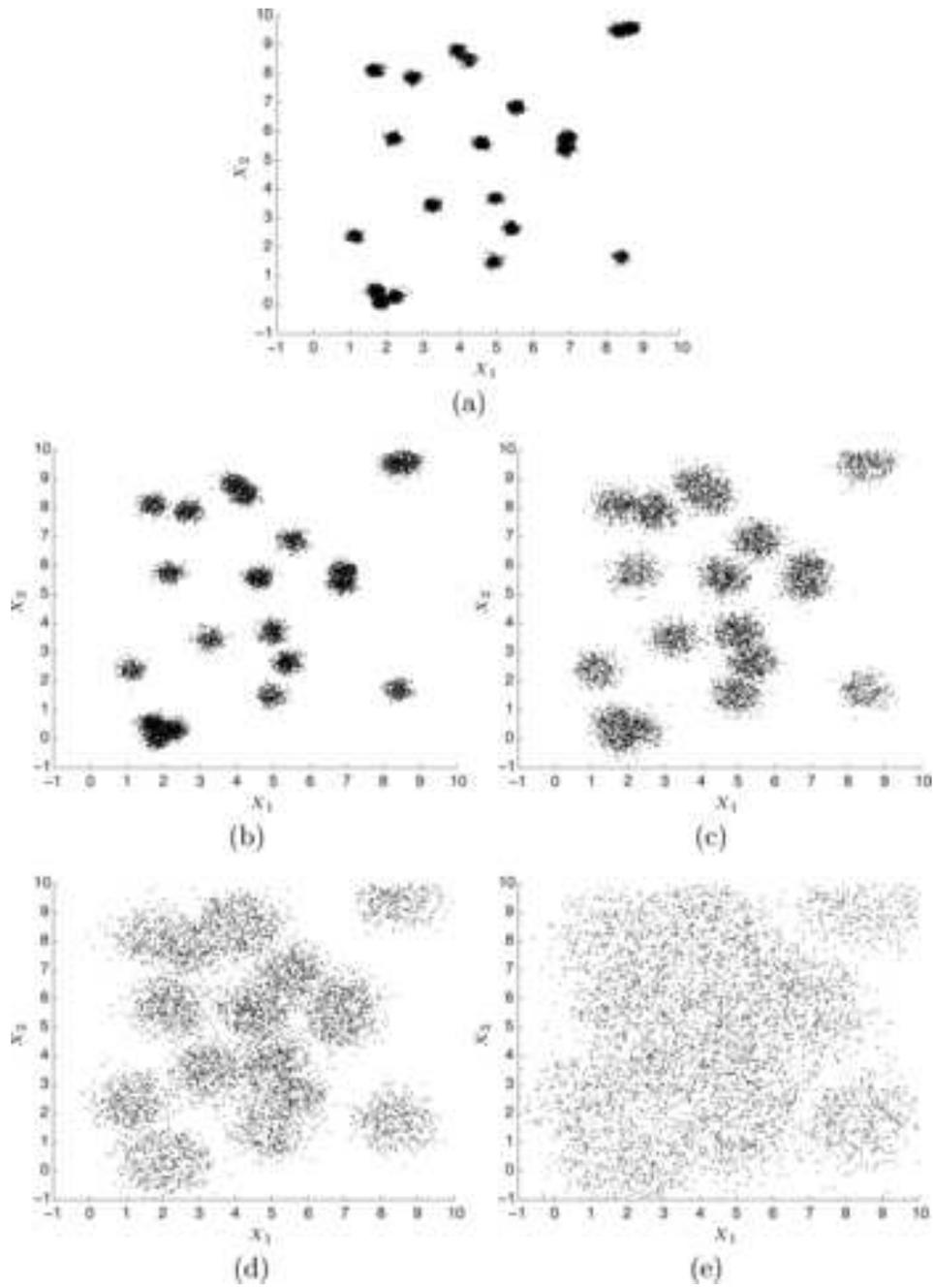

FIG. 2.   *Samples generated from each chain of the EE sampler. The chains are sorted in ascending order of temperature and energy truncation.*

To compare the mixing speed for the two sampling algorithms, we counted in each of the 20 runs how many times the samples visited each mode in the last 2000 iterations. We then calculated the absolute frequency error for each mode, $\mathrm{err}_i = |\hat{f}_i - 0.05|$, where $\hat{f}_i$ is the sample frequency of the $i$th mode $(i = 1, \ldots, 20)$ being visited. For each mode $i$, we calculated the median and the maximum of $\mathrm{err}_i$ over the 20 runs. Table 3 reports, for each mode, the ratio $R1$ of the median frequency error of PT over that of EE, and the ratio $R2$ of the maximum frequency error of PT over the corresponding value of EE. A two- to fourfold improvement by EE was observed. We also noted that EE did not miss a single mode in all runs, whereas PT missed some modes in each run. Table 3 also shows, out of the 20 runs, how many times each mode was missed by PT; for example, mode 1 was missed twice by PT over the 20 runs. The better global exploration ability of the EE sampler is thus evident. To further compare the two algorithms by their convergence speed, we tuned the temperature ladder for PT to achieve the *best* performance of PT; for example, we tuned the highest temperature $T_4$ of PT in the range $[5, 100]$. We observe that the best performance of PT such that it is not trapped by some local modes in 50,000 samples and that it achieves minimal sample autocorrelation is the setting associated with $T_4 = 10$. The sample autocorrelations of the target chain $X^{(0)}$ from the EE sampler and the optimal PT are shown in Figures 3(c) and (d), respectively; evidently the autocorrelation of the EE sampler decays much faster even compared with a well-tuned PT.

We also use this example to study the choice of $p_{\mathrm{ee}}$. We took $p_{\mathrm{ee}} = 0.1$, 0.2, 0.3 and 0.4, and ran the EE sampler 20 times independently for each value of $p_{\mathrm{ee}}$. From the average mean squared errors for estimating $(EX_1, EX_2)$ and $(EX_1^2, EX_2^2)$ we find that the EE sampler behaves well when $p_{\mathrm{ee}}$ is between 0.1 and 0.3. When $p_{\mathrm{ee}}$ is increased to 0.4, the performance of the EE sampler worsens. In addition, we also noticed that the performance of the EE sampler is not sensitive to the value of $p_{\mathrm{ee}}$, as long as it is in the range $[0.05, 0.3]$.

Next, we changed the weight and variance for each component in (6) such that $w_i \propto 1/d_i$ and $\sigma_i^2 = d_i/20$, where $d_i = \|\boldsymbol{\mu}_i - (5, 5)'\|$. The distributions closer to $(5, 5)$ have larger weights and lower energy. We used this example to test our strategy of dynamically updating the energy and temperature ladders (Section 3.3). We set the initial energy lower bound $H_0 = 3$, which is higher than the energy at any of the 20 modes (in practice, we could try to get a better initial value of $H_0$ by some local optimization). The

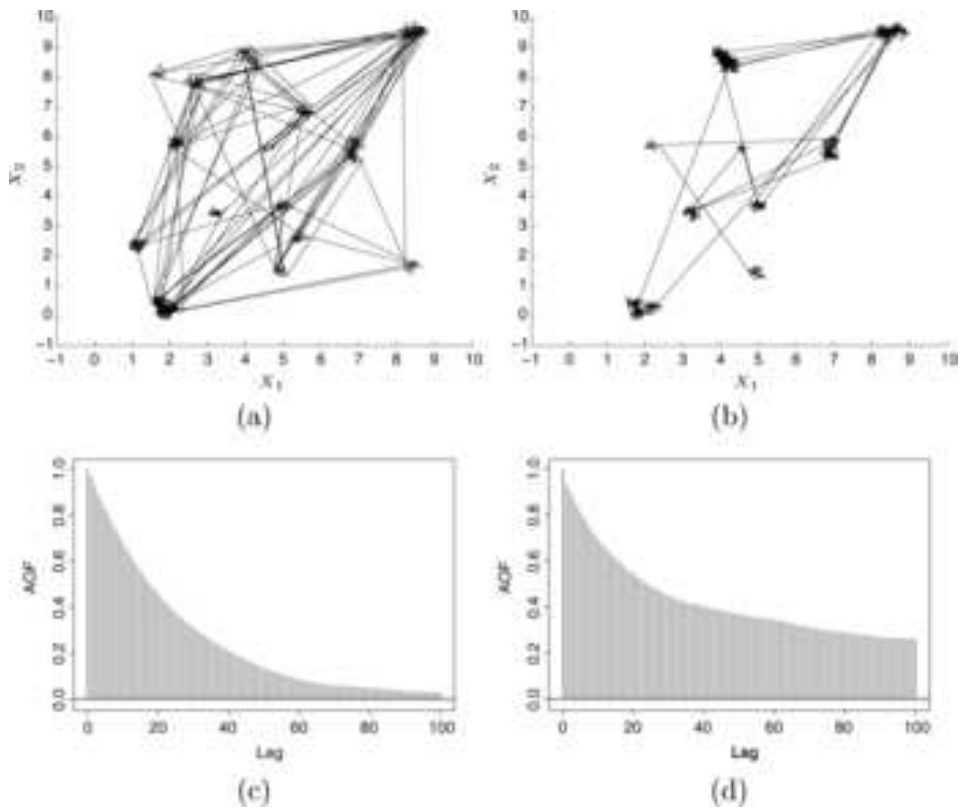

Fig. 3. *Mixture normal distribution with equal weights and variances. The sample path of the last* 2000 *iterations for* (a) *EE sampler and* (b) *parallel tempering. Autocorrelation plots for the samples from* (c) *EE sampler and* (d) *optimal parallel tempering.*

TABLE 2
*Comparison of the EE sampler and PT for estimating the mixture normal distributions*

|  | $\boldsymbol{EX_1}$ | $\boldsymbol{EX_2}$ | $\boldsymbol{EX_1^2}$ | $\boldsymbol{EX_2^2}$ |
|---|---|---|---|---|
| True value | 4.478 | 4.905 | 25.605 | 33.920 |
| EE | 4.5019 (0.107) | 4.9439 (0.139) | 25.9241 (1.098) | 34.4763 (1.373) |
| PT | 4.4185 (0.170) | 4.8790 (0.283) | 24.9856 (1.713) | 33.5966 (2.867) |
| MSE(PT)/MSE(EE) | 2.7 | 3.8 | 2.6 | 3.8 |
| True value | 4.688 | 5.030 | 25.558 | 31.378 |
| EE | 4.699 (0.072) | 5.037 (0.086) | 25.693 (0.739) | 31.433 (0.839) |
| PT | 4.709 (0.116) | 5.001 (0.134) | 25.813 (1.122) | 31.105 (1.186) |
| MSE(PT)/MSE(EE) | 2.6 | 2.5 | 2.4 | 2.1 |

The numbers in parentheses are the standard deviations from 20 independent runs. The upper and bottom halves correspond to equal and unequal weights and variances, respectively.

Table 3
*Comparison of mixing speed of the EE sampler and PT for the mixture normal distribution*

|  | $\mu_1$ | $\mu_2$ | $\mu_3$ | $\mu_4$ | $\mu_5$ | $\mu_6$ | $\mu_7$ | $\mu_8$ | $\mu_9$ | $\mu_{10}$ |
|---|---|---|---|---|---|---|---|---|---|---|
| $R1$ | 2.8 | 2.9 | 3.2 | 3.6 | 2.7 | 4.1 | 1.7 | 2.1 | 2.7 | 3.3 |
| $R2$ | 1.6 | 4.5 | 2.5 | 4.9 | 2.5 | 2.1 | 2.2 | 2.3 | 2.3 | 2.5 |
| $PT_{\mathrm{mis}}$ | 2 | 3 | 1 | 5 | 4 | 0 | 2 | 2 | 3 | 1 |

|  | $\mu_{11}$ | $\mu_{12}$ | $\mu_{13}$ | $\mu_{14}$ | $\mu_{15}$ | $\mu_{16}$ | $\mu_{17}$ | $\mu_{18}$ | $\mu_{19}$ | $\mu_{20}$ |
|---|---|---|---|---|---|---|---|---|---|---|
| $R1$ | 1.3 | 2.5 | 2.0 | 4.4 | 2.2 | 2.5 | 2.0 | 1.9 | 2.6 | 1.4 |
| $R2$ | 3.8 | 1.2 | 1.9 | 2.1 | 3.2 | 3.0 | 3.1 | 3.0 | 3.9 | 2.0 |
| $PT_{\mathrm{mis}}$ | 1 | 3 | 1 | 1 | 1 | 3 | 4 | 6 | 1 | 3 |

$R1$ is the ratio of median frequency error of PT over that of EE; $R2$ is the ratio of maximum frequency error of PT over that of EE. $PT_{\mathrm{mis}}$ reports the number of runs for which PT missed the individual modes. The EE sampler did not miss any of the 20 modes. All the statistics are calculated from the last 2000 samples in 20 independent runs.

highest energy level and temperature were set at 100 and 20, respectively. We started with five chains and dynamically added more chains if necessary. The strategy for automatically updating the proposal step size was applied as well. After drawing the samples we also calculated the first two sample moments from the target chain $X^{(0)}$ as simple estimates for the theoretical moments. The mean and standard deviation of the estimates based on 20 independent EE runs, each consisting of 10,000 iterations, are shown in Table 2 (the bottom half ). The sample path for the last 1000 iterations and the autocorrelation plots are shown in Figure 4.

For comparison, PT was applied to this unequal weight case as well. PT used the same total number of chains as the EE sampler. The MH step size of PT was tuned to achieve the same acceptance rate. The temperature ladder of PT was also tuned so that the exchange operator in PT had the same acceptance rate as the equi-energy jump in EE, similarly to what we did in the previous comparison. With these well-tuned parameters, we ran PT for the same number of iterations and calculated the first two sample moments as we did for the EE samples. The results are reported in Table 2. It is seen that the EE sampler with the self-adaptation strategies provided more precise estimates (both smaller bias and smaller variance) in all the cases. Similar improvements in mixing speed and autocorrelation decay were also observed (Figure 4).

**4. Calculating the density of states and the Boltzmann averages.** The previous section illustrated the benefit of constructing the energy rings: It allows more efficient sampling through the equi-energy jump. By enabling one to look at the states within a given energy range, the energy rings also provide a direct means to study the microscopic structure of the state space on an energy-by-energy basis, that is, the microcanonical distribution. The EE sampler and the energy rings are thus well suited to study problems in statistical mechanics.

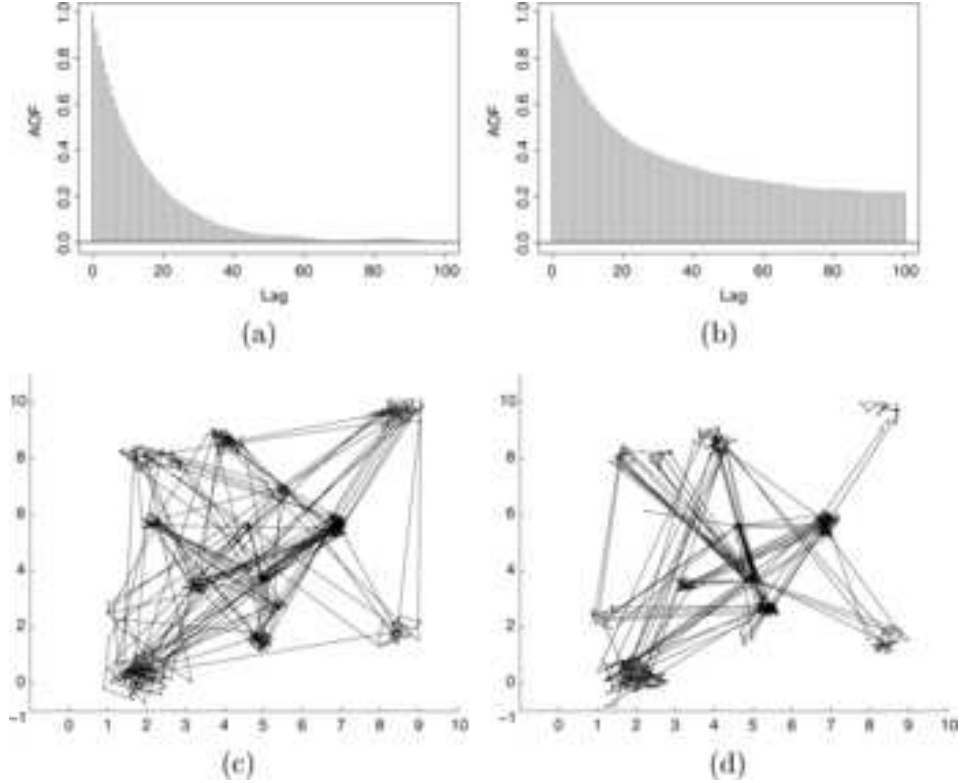

Fig. 4.  *Mixture normal distribution with unequal weights and variances. Autocorrelation plots for the samples from (*a*) EE sampler and (*b*) parallel tempering. The sample path for the last* 1000 *iterations for (*c*) EE sampler and (*d*) parallel tempering.*

Starting from the Boltzmann distribution (1) one wants to study various aspects of the system. The density of states $\Omega(u)$, whose logarithm is referred to as the microcanonical entropy, plays an important role in the study, because in addition to the temperature–energy duality depicted in Section 2, many thermodynamic quantities can be directly calculated from the density of states, for example, the partition function $Z(T)$, the internal energy $U(T)$, the specific heat $C(T)$ and the free energy $F(T)$ can be calculated via

$$Z(T) = \int \Omega(u) e^{-u/T} \, du,$$

$$U(T) = \frac{\int u\,\Omega(u)e^{-u/T}\,du}{\int \Omega(u)e^{-u/T}\,du},$$

$$C(T) = \frac{\partial U(T)}{\partial T} = \frac{1}{T^2}\left[\frac{\int u^2\Omega(u)e^{-u/T}\,du}{\int \Omega(u)e^{-u/T}\,du} - \left(\frac{\int u\Omega(u)e^{-u/T}\,du}{\int \Omega(u)e^{-u/T}\,du}\right)^2\right],$$

$$F(T) = -T\log(Z(T)).$$

Since the construction of the energy rings is an integral part of the EE sampler, it leads to a simple way to estimate the density of states. Suppose we have a discrete system to study, and after performing the EE sampling on the distributions

$$\pi_i(x) \propto \exp(-h_i(x)), \qquad h_i(x) = \frac{1}{T_i}(h(x) \vee H_i), \qquad 0 \le i \le K,$$

we obtain the energy rings $\hat{D}_j^{(i)}$ ($0 \le i, j \le K$), each $\hat{D}_j^{(i)}$ corresponding to an a priori energy range $[H_j, H_{j+1}]$. To calculate the density of states for the discrete system we can further divide the energy rings into subsets such that each subset corresponds to one energy level. Let $m_{iu}$ denote the total number of samples in the $i$th chain $X^{(i)}$ that have energy $u$. Clearly $\sum_{u \in [H_j, H_{j+1})} m_{iu} = |\hat{D}_j^{(i)}|$. Under the distribution $\pi_i$

$$(7) \qquad P_{\pi_i}(h(X) = u) = \frac{\Omega(u)e^{-(u \vee H_i)/T_i}}{\sum_v \Omega(v)e^{-(v \vee H_i)/T_i}}.$$

Since the density of states $\Omega(u)$ is common for each $\pi_i$, we can combine the sample chains $X^{(i)}$, $0 \le i \le K$, to estimate $\Omega(u)$. Denote $m_{i\bullet} = \sum_u m_{iu}$, $m_{\bullet u} = \sum_i m_{iu}$ and $a_{iu} = e^{-(u \vee H_i)/T_i}$ for notational ease. Pretending we have independent multinomial observations,

$$(\ldots, m_{iu}, \ldots) \sim \text{multinomial}\left(m_{i\bullet}; \ldots, \frac{\Omega(u)a_{iu}}{\sum_v \Omega(v)a_{iv}}, \ldots\right),$$

$$(8)$$

$$i = 0, 1, \ldots, K,$$

the MLE of $\Omega(u)$ is then

$$(9) \qquad \hat{\Omega} = \arg\max_{\Omega}\left\{\sum_u m_{\bullet u}\log(\Omega(u)) - \sum_i m_{i\bullet}\log\left(\sum_v \Omega(v)a_{iv}\right)\right\}.$$

Since $\Omega(u)$ is specified up to a scale change [see (7)], to estimate the relative value we can without loss of generality set $\Omega(u_0) = 1$ for some $u_0$. The first-order condition of (9) gives

$$(10) \qquad \frac{m_{\bullet u}}{\hat{\Omega}(u)} - \sum_i \frac{m_{i\bullet}a_{iu}}{\sum_v \hat{\Omega}(v)a_{iv}} = 0 \qquad \text{for all } u,$$

which can be used to compute $\hat{\Omega}(u)$ through a simple iteration,

$$(11) \qquad \hat{\Omega}(u) = m_{\bullet u} \Big/ \sum_i \frac{m_{i\bullet} a_{iu}}{\sum_v \hat{\Omega}(v) a_{iv}}.$$

A careful reader might question the independent multinomial assumption (9). But it is only used to motivate (10), which itself can be viewed as a moment equation and is valid irrespective of the multinomial assumption.

With the density of states estimated, suppose one wants to investigate how the Boltzmann average $E(g(X);T) = \frac{\sum_x g(x)\exp(-h(x)/T)}{\sum_x \exp(-h(x)/T)}$ varies as a function of temperature $T$ (e.g., phase transition). Then we can write (see Lemma 1)

$$E(g(X);T) = \frac{\sum_u \Omega(u)e^{-u/T}\nu_g(u)}{\sum_u \Omega(u)e^{-u/T}}.$$

To estimate the microcanonical average $\nu_g(u) = E(g(X)|h(X) = u)$, we can simply calculate the sample average over the energy slice $\{x : h(x) = u\}$ for each chain $X^{(0)}, \ldots, X^{(K)}$ and combine them using weights proportional to the energy slice sample size $m_{iu}$ for $i = 0, 1, \ldots, K$.

So far we have focused on discrete systems. In continuous systems to estimate $\Omega(u)$ and $\nu_g(u) = E(g(X)|h(X) = u)$ we can simply discretize the energy space with an acceptable resolution and follow the preceding approach to use the EE sampler and the energy rings.

To illustrate our method to calculate the density of states and the Boltzmann averages, let us consider a multidimensional normal distribution $f(\mathbf{x};T) = \frac{1}{Z(T)}\exp[-h(\mathbf{x})/T]$ with temperature $T$ and the energy function $h(\mathbf{x}) = \frac{1}{2}\sum_i^n x_i^2$. This corresponds to a system of $n$ uncoupled harmonic oscillators. Since the analytical result is known in this case, we can check our numerical calculation with the exact values. Suppose we are interested in estimating $E(X_1^2;T)$ and the ratio of normalizing constants $Z(T)/Z(1)$ for $T$ in $[1, 5]$. The theoretical density of states (from the $\chi^2$ result in this case) is $\Omega(u) \propto u^{n/2-1}$ and the microcanonical average $\nu_g(u) = E(X_1^2|h(\mathbf{X}) = u) = 2u/n$. Our goal is to estimate $\Omega(u)$ and $\nu_g(u)$ and then combine them to estimate both $E(X_1^2;T)$ and $Z(T)/Z(1)$ as functions of the temperature $T$.

We took $n = 4$ and applied the EE sampler with five chains to sample from the four-dimensional distribution. The energy levels $H_0, H_1, \ldots, H_4$ were assigned by a geometric progression along the interval $[0, 50]$. The temperatures were accordingly set between 1 and 20. The equi-energy jump probability $p_{ee}$ was taken to be 0.05. Each chain was run for 150,000 iterations (the first 50,000 iterations were the burn-in period); the sampling was repeated ten times to calculate the standard error of our estimates. Since the underlying distribution is continuous, to estimate $\Omega(u)$ and $\nu_g(u)$ each energy ring was further evenly divided into 20 small intervals. The estimates $\hat{\Omega}(u)$ and $\hat{\nu}_g(u)$ were calculated from the recursion (11) and the combined energy slice sample average, respectively. Figure 5(a) and (b) shows $\hat{\Omega}(u)$ and $\hat{\nu}_g(u)$ compared with their theoretical values. Using $\hat{\Omega}(u)$ and $\hat{\nu}_g(u)$, we estimated $E(X_1^2; T)$ and $Z(T)/Z(1)$ by

$$E(X_1^2; T) \approx \frac{\sum_u \hat{\nu}_g(u)\hat{\Omega}(u)e^{-u/T}\Delta u}{\sum_u \hat{\Omega}(u)e^{-u/T}\Delta u},$$

$$\frac{Z(T)}{Z(1)} \approx \frac{\sum_u \hat{\Omega}(u)e^{-u/T}\Delta u}{\sum_u \hat{\Omega}(u)e^{-u}\Delta u}.$$

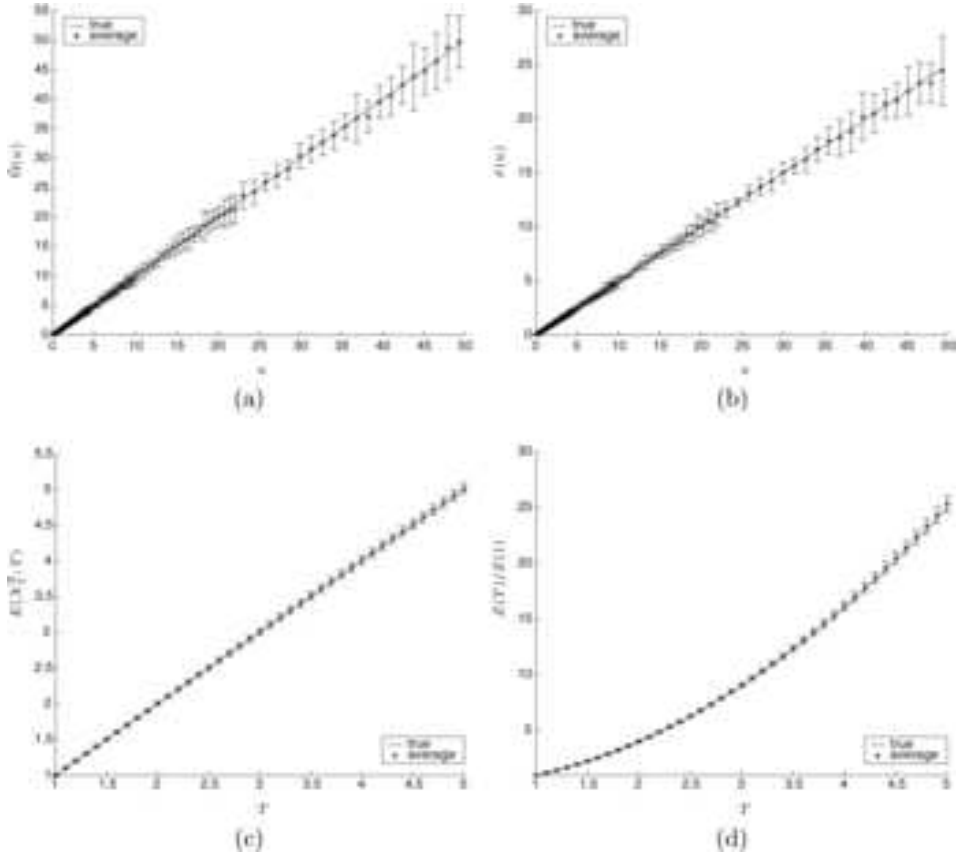

FIG. 5. *The density of states calculation for the four-dimensional normal distribution. (a) Estimated density of states $\hat{\Omega}(u)$; (b) estimated micro-canonical average $\hat{\nu}(u)$; (c) estimated $E(X_1^2; T)$ and (d) estimated $Z(T)/Z(1)$ for $T$ varying from 1 to 5. The error bars represent plus or minus twice the standard deviation from ten independent runs.*

Letting $T$ vary in $[1, 5]$ results in the estimated curves in Figure 5(c) and (d). One can see from the figure that our estimations are very precise compared to the theoretical values. In addition, our method has the advantage that we are able to construct estimates for a wide temperature range using one simulation that involves only five temperatures.

We further tested our method on a four-dimensional normal mixture distribution with the energy function

$$(12) \qquad h(\mathbf{x}) = -\log[\exp(-\|\mathbf{x} - \boldsymbol{\mu}_1\|^2) + 0.25 \exp(-\|\mathbf{x} - \boldsymbol{\mu}_2\|^2)],$$

where $\boldsymbol{\mu}_1 = (3, 0, 0, 0)'$ and $\boldsymbol{\mu}_2 = (-3, 0, 0, 0)'$. For the Boltzmann distribution $\frac{1}{Z(T)} \exp(-h(\mathbf{x})/T)$, we are interested in estimating the probability $P(X_1 > 0; T)$ and studying how it varies as $T$ changes. It is easy to see if $T = 1$ this probability will be $0.8 = 1/1.25$, the mixture proportion, and it decreases as $T$ becomes larger. We applied the EE sampler to this problem with the same parameter setting as in the preceding 4D normal example. Using the energy rings, we calculated the estimates $\hat{\Omega}(u)$ and $\hat{\nu}_g(u)$ and then combined them to estimate $P(X_1 > 0; T)$. Figure 6 plots the estimates. It is interesting to note from the figure that the density of states for the mixture model has a change point at energy $u = 1.4$. This is due to the fact that the energy at the mode $\boldsymbol{\mu}_2$ is about $1.4 (\approx -\log 0.25)$, and hence for $u < 1.4$ all the samples are from the first mode $\boldsymbol{\mu}_1$, and for $u > 1.4$ the samples come from both modes, whence a change point appears. A similar phenomenon occurs in the plot for $\nu_g(u) = P(X_1 > 0 | h(\mathbf{X}) = u)$. The combined estimate of $P(X_1 > 0; T)$ in Figure 6 was checked and agreed very well with the exact values from numerical integration.

**5. Statistical estimation with the EE sampler.** In statistical inference, usually after obtaining the Monte Carlo samples the next goal is to estimate some statistical quantities. Unlike statistical mechanical considerations, in statistical inference problems often one is only interested in one target distribution (i.e., only one temperature). Suppose the expected value $E_{\pi_0} g(X)$ under the target distribution $\pi_0 = \pi$ is of interest. A simple estimate is the sample mean under the chain $X^{(0)}$ ($T = 1$). This, however, does not use the samples optimally in that it essentially throws away all the other sampling chains $X^{(1)}, \ldots, X^{(K)}$. With the help of the energy rings, in fact all the chains can be combined together to provide a more efficient estimate. One way is to use the energy-temperature duality as discussed above. Here we present an alternative, more direct method: We directly work with the (finite number of ) expectations within each energy set. For a continuous problem, this allows us to avoid dealing with infinitesimal energy intervals, which are needed in the calculation of microcanonical averages.

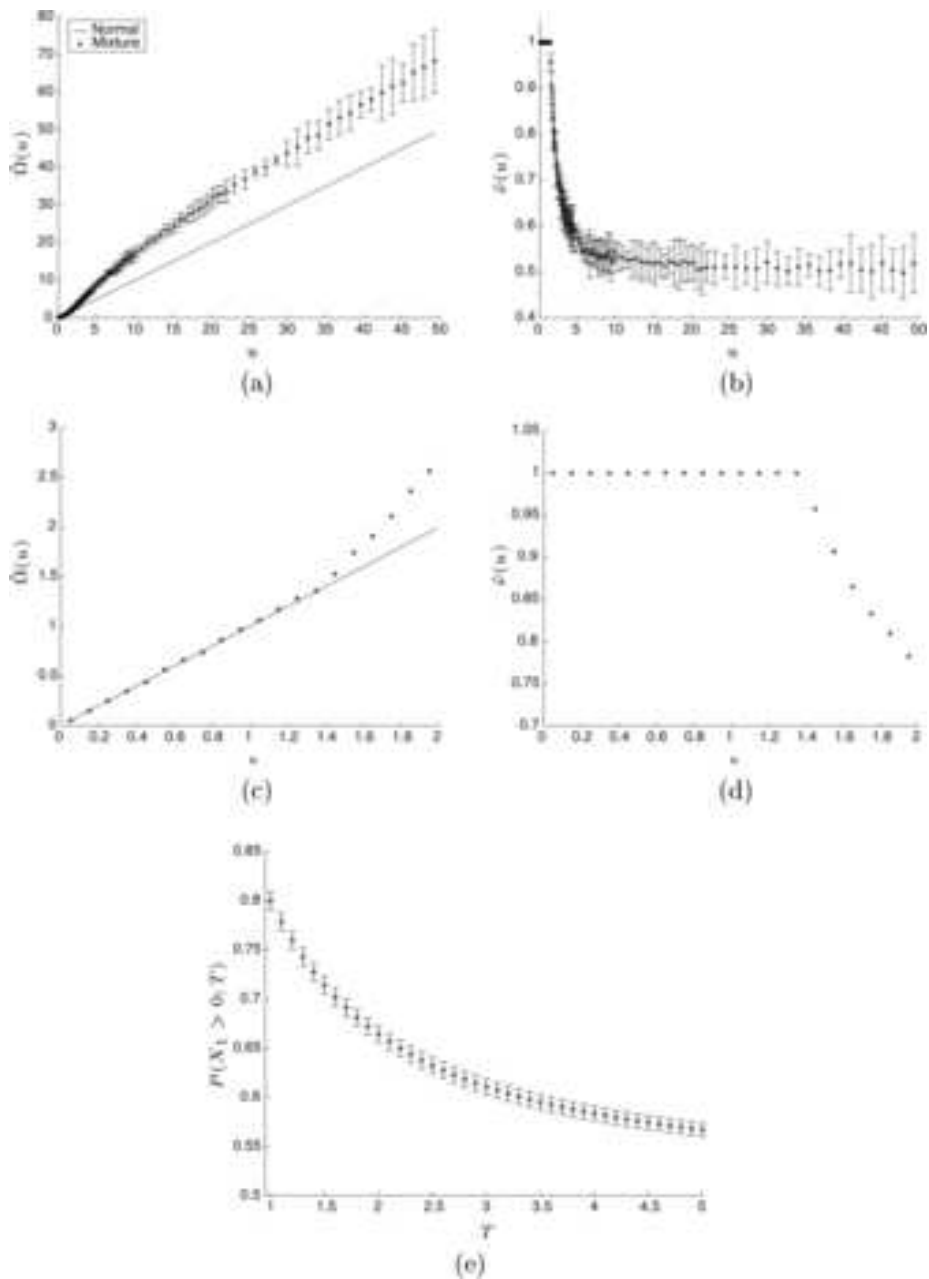

Fig. 6.  *Four-dimensional normal mixture distribution. (a) Estimated density of states $\hat{\Omega}(u)$; (b) estimated microcanonical average $\hat{\nu}(u)$; (c) detailed plot for $\hat{\Omega}(u)$ around the change point; (d) detailed plot for $\hat{\nu}(u)$ around the change point; (e) estimated $P(X_1 > 0; T)$ for $T$ in $[1, 5]$. For comparison, the theoretical $\Omega(u)$ from the 4D normal example is also shown in (a) and (c). The error bars represent plus or minus twice the standard deviation from ten independent runs.*

The starting point is the identity

$$(13) \qquad E_{\pi_0} g(X) = \sum_{j=0}^{K} p_j E_{\pi_0} \{g(X) | X \in D_j\},$$

where $p_j = P_{\pi_0}(X \in D_j)$, which suggests that we can first estimate $p_j$ and $G_j = E_{\pi_0} \{g(X) | X \in D_j\}$, and then conjoin them. A naive estimate of $G_j$ is the sample average of $g(x)$ within the energy ring $\hat{D}_j^{(0)}$. But a better way is to use the energy rings $\hat{D}_j^{(i)}$ for $i = 0, 1, \dots, K$ together, because for each $i$ the importance-weighted

$$\hat{G}_j^{(i)} = \frac{\sum_{X \in \hat{D}_j^{(i)}} g(X) w^{(i)}(X)}{\sum_{X \in \hat{D}_j^{(i)}} w^{(i)}(X)},$$

where $w^{(i)}(x) = \exp\{h_i(x) - h(x)\}$, is a consistent estimate of $G_j$. We can thus put proper weights on them to combine them.

Since it is quite often that a number of different estimands $g$ are of interest, a conceptually simple and effective method is to weight $\hat{G}_j^{(i)}$, $i = 0, 1, \dots, K$, proportionally to their *effective sample sizes* [15],

$$ESS_j^{(i)} = \frac{|\hat{D}_j^{(i)}|}{1 + \mathrm{Var}_{\pi_i}\{w^{(i)}(X) | X \in D_j\} / (E_{\pi_i}\{w^{(i)}(X) | X \in D_j\})^2},$$

which measures the effectiveness of a given importance-sampling scheme, and does not involve the specific estimand. $ESS_j^{(i)}$ can be easily estimated by using the sample mean and variance of $w^{(i)}(X)$ for $X$ in the energy ring $\hat{D}_j^{(i)}$.

To estimate the energy-ring probability $p_j = P_{\pi_0}(X \in D_j)$, one can use the sample size proportion $\frac{|\hat{D}_j^{(0)}|}{\sum_{k=0}^{K} |\hat{D}_k^{(0)}|}$. But again it is better to use all the chains, because for any $i$,

$$\hat{p}_j^{(i)} = \frac{\sum_{X \in \hat{D}_j^{(i)}} w^{(i)}(X)}{\sum_{X \in \hat{D}^{(i)}} w^{(i)}(X)}, \qquad \text{where } \hat{D}^{(i)} = \bigcup_j \hat{D}_j^{(i)}$$

is a consistent estimate of $p_j$. To combine them properly, we can weight them inversely proportional to their variances. A simple delta method calculation (disregarding the dependence) gives the asymptotic variance of $\hat{p}_j^{(i)}$,

$$V_j^{(i)} = \frac{1}{n_i} \frac{E_{\pi_i}[(I(X \in D_j) - p_j) w^{(i)}(X)]^2}{[E_{\pi_i} w^{(i)}(X)]^2},$$

where $n_i$ is the total number of samples under the $i$th chain $X^{(i)}$. $V_j^{(i)}$ can be estimated by

$$\hat{V}_j^{(i)} = \frac{\sum_{X \in \hat{D}^{(i)}} [(I(X \in D_j) - \tilde{p}_j) w^{(i)}(X)]^2}{(\sum_{X \in \hat{D}^{(i)}} w^{(i)}(X))^2}$$

$$= (1 - 2\tilde{p}_j) \frac{\sum_{X \in \hat{D}_j^{(i)}} [w^{(i)}(X)]^2}{(\sum_{X \in \hat{D}^{(i)}} w^{(i)}(X))^2} + \tilde{p}_j^2 \frac{\sum_{X \in \hat{D}^{(i)}} [w^{(i)}(X)]^2}{(\sum_{X \in \hat{D}^{(i)}} w^{(i)}(X))^2},$$

where $\tilde{p}_j$ is a consistent estimate of $p_j$. Since $\hat{V}_j^{(i)}$ requires an a priori $\tilde{p}_j$, a simple iterative procedure can be applied to obtain a good estimate of $p_j$, starting from the naive $\hat{p}_j^{(0)}$. In practice to ensure the numerical stability of the variance estimate $\hat{V}_j^{(i)}$ it is recommended that $\hat{p}_j^{(i)}$ be included in the combined estimate only if the sample size $|\hat{D}_j^{(i)}|$ is reasonably large, say more than 50. When we get our estimates for $p_j$, we need to normalize them before plugging in (13) to calculate our combined estimates.

Since our energy-ring estimate of $E_{\pi_0} g(X)$ employs all the sampling chains, one expects it to increase the estimation efficiency by a large margin compared with the naive method.

*The two-dimensional normal mixture model* (*continued*). To illustrate our estimation strategy, consider again the 2D mixture model (6) in Section 3.4 with equal weights and variances $\sigma_1 = \cdots = \sigma_{20} = 0.1$, $w_1 = \cdots = w_{20} = 0.05$.

Suppose we want to estimate the functions $EX_1^2$, $EX_2^2$, $Ee^{-10X_1}$ and $Ee^{-10X_2}$ and the tail probabilities

$$p_1 = P(X_1 > 8.41, X_2 < 1.68 \text{ and } \sqrt{(X_1 - 8.41)^2 + (X_2 - 1.68)^2} > 4\sigma),$$

TABLE 4
*Comparison of the energy-ring estimates with the naive estimates that use $X^{(0)}$ only*

|  | $EX_1^2$ | $EX_2^2$ | $Ee^{-10X_1}$ | $Ee^{-10X_2}$ | $p_1$ | $p_2$ |
|---|---|---|---|---|---|---|
| **True value** | **25.605** | **33.920** | **9.3e-7** | **0.0378** | **4.2e-6** | **6.7e-5** |
| Energy ring estimates | 25.8968 (0.9153) | 34.2902 (1.1579) | 8.8e-7 (1.2e-7) | 0.0379 (0.0044) | 4.5e-6 (1.5e-6) | 6.4e-5 (2.0e-5) |
| Naive estimates | 25.9241 (1.0982) | 34.4763 (1.3733) | 8.7e-7 (1.5e-7) | 0.0380 (0.0052) | 1.0e-5 (2.5e-5) | 7.3e-5 (6.2e-5) |
| MSER | 71% | 67% | 57% | 72% | 0.34% | 11% |

The numbers in parentheses are the standard deviations from 20 independent runs. MSER is defined as $MSE_1/MSE_2$, where $MSE_1$ and $MSE_2$ are the mean squared errors of the energy-ring estimates and the naive estimates, respectively.

$p_2 = P(X_1^2 + X_2^2 > 175)$.

After obtaining the EE samples (as described in Section 3.4), we calculated our energy-ring estimates. For comparison, we also calculated the naive estimates based on the target chain $X^{(0)}$ only. The calculation was repeated for the 20 independent runs of the EE sampler under the same parameter settings. Table 4 reports the result. Evidently, the energy-ring estimates with both smaller bias *and* smaller variance are much more precise than the naive estimates. The mean squared error has been reduced by at least 28% in all cases. The improvement over the naive method is particularly dramatic in the tail probability estimation, where the MSEs experienced a more than ninefold reduction.

**6. Applications of the EE sampler.** We will apply the EE sampler and the estimation strategy to a variety of problems in this section to illustrate its effectiveness. The first example is a mixture regression problem. The second example involves motif discovery in computational biology. In the third one we study the thermodynamic property of a protein folding model.

6.1. *Mixture exponential regression.* Suppose we observe data pairs $(\mathbf{Y}, \mathbf{X}) = \{(y_1, \mathbf{x}_1), \ldots, (y_n, \mathbf{x}_n)\}$ from the mixture model

$$(14) \qquad y_i \sim \begin{cases} \text{Exp}[\theta_1(\mathbf{x}_i)], & \text{with probability } \alpha, \\ \text{Exp}[\theta_2(\mathbf{x}_i)], & \text{with probability } 1 - \alpha, \end{cases}$$

where $\theta_j(\mathbf{x}_i) = \exp[\boldsymbol{\beta}_j^T \mathbf{x}_i]$ $(j = 1, 2)$, $\text{Exp}(\theta)$ denotes an exponential distribution with mean $\theta$ and $\boldsymbol{\beta}_1$, $\boldsymbol{\beta}_2$ and $\mathbf{x}_i$ $(i = 1, 2, \ldots, n)$ are $p$-dimensional vectors. Given the covariates $\mathbf{x}_i$ and the response variable $y_i$, one wants to infer the regression coefficients $\boldsymbol{\beta}_1$ and $\boldsymbol{\beta}_2$ and the mixture probability $\alpha$. The likelihood of the observed data is

$$(15) \qquad \begin{aligned} &P(\mathbf{Y}, \mathbf{X} | \alpha, \boldsymbol{\beta}_1, \boldsymbol{\beta}_2) \\ &\quad \propto \prod_{i=1}^{n} \left[ \frac{\alpha}{\theta_1(\mathbf{x}_i)} \exp\left(-\frac{y_i}{\theta_1(\mathbf{x}_i)}\right) + \frac{1 - \alpha}{\theta_2(\mathbf{x}_i)} \exp\left(-\frac{y_i}{\theta_2(\mathbf{x}_i)}\right) \right]. \end{aligned}$$

If we put a Beta$(1, 1)$ prior on $\alpha$, and a multivariate normal $N(\mathbf{0}, \sigma^2 \mathbf{I})$ on $\boldsymbol{\beta}_j$ $(j = 1, 2)$, the energy function, defined as the negative log-posterior density, is

$$(16) \qquad \begin{aligned} h(\alpha, \boldsymbol{\beta}_1, \boldsymbol{\beta}_2) &= -\log P(\alpha, \boldsymbol{\beta}_1, \boldsymbol{\beta}_2 | \mathbf{Y}, \mathbf{X}) \\ &= -\log P(\mathbf{Y}, \mathbf{X} | \alpha, \boldsymbol{\beta}_1, \boldsymbol{\beta}_2) + \frac{1}{2\sigma^2} \sum_{k=1}^{2} \sum_{j=1}^{p} \boldsymbol{\beta}_{kj}^2 + C, \end{aligned}$$

where $C$ is a constant. Since $h(\alpha, \boldsymbol{\beta}_1, \boldsymbol{\beta}_2) = h(1 - \alpha, \boldsymbol{\beta}_2, \boldsymbol{\beta}_1)$ (i.e., nonidentifiable), the posterior distribution has multiple modes in the $(2p + 1)$-dimensional parameter space. This example thus serves as a good model to test the performance of the EE sampler in high-dimensional multimodal problems.

We simulated 200 data pairs with the following parameter setting: $\alpha = 0.3$, $\boldsymbol{\beta}_1 = (1, 2)'$, $\boldsymbol{\beta}_2 = (4, 5)'$ and $\mathbf{x}_i = (1, u_i)'$, with the $u_i$'s independently drawn from Unif$(0, 2)$. We took $\sigma^2 = 100$ in the prior normal distributions for the regression coefficients. After a local minimization of the energy function (16) from several randomly chosen states in the parameter space, the minimum energy $H_{\min}$ is found to be around $H_{\min} \approx -1740.8$ (this value is not crucial in the EE sampler, since it can adaptively adjust the energy and temperature ladder; see Section 3.3). We then applied the EE sampler with eight chains to sample from the posterior distribution. The energy ladder was set between $H_{\min}$ and $H_{\min} + 100$ in a geometric progression, and the temperatures were between 1 and 30. The equi-energy jump probability $p_{\text{ee}}$ was taken to be 0.1. Each chain was run for 15,000 iterations with the first 5000 iterations serving as a burn-in period. The overall acceptance rates for the MH move and the equi-energy jump were 0.27 and 0.80, respectively. Figure 7(a) to (c) shows the sample marginal posterior distributions of $\alpha$, $\boldsymbol{\beta}_1$ and $\boldsymbol{\beta}_2$ from the target chain $X^{(0)}$. It is clear that the EE sampler visited the two modes, equally frequently in spite of the long distance between the two modes, and the samples around each mode were centered at the true parameters. Furthermore, since the posterior distribution for this problem has two symmetric modes in the parameter space due to the nonidentifiability, by visiting the two modes equally frequently the EE sampler demonstrates its capability of global exploration (as opposed to being trapped by a local mode).

For comparison, we also applied PT with eight chains to this problem under the same parameter setting, where the acceptance rates for the MH move and the exchange operator were 0.23 and 0.53, respectively. To compare their ability to escape local traps, we calculated the frequency with which the samples stayed at one particular mode in the lowest temperature chain ($T_0 = 1$). This frequency was 0.55 for the EE samplers and 0.89 for PT, indicating that the EE sampler visited the two modes symmetrically, while PT tended to be trapped at one mode for a long time. We further tuned the temperature ladder for PT with the highest temperature varying from 10 to 50. For each value of the highest temperature, the temperature ladder was set to decrease with a geometric rate. We observed that the sample autocorrelations decrease with the decrease of the highest temperature, that is, the denser the temperature ladder, the smaller the sample autocorrelation. However, the tradeoff is that with lower temperature, PT tends to get trapped in one local mode. For instance, PT was totally trapped to one particular mode if we set the highest temperature to be 10 although with this temperature ladder PT showed autocorrelation comparable with that of EE. On the other hand, if we increased the temperatures, PT would be able to jump between the two modes, but the autocorrelation also increased since the exchange rates became lower. But even with the highest temperature raised to 50, 80% of the PT samples were trapped to one mode and the posterior distributions were severely asymmetric. The autocorrelation plots for $\beta_{11}$ are compared in Figure 7(d) and (e), where one sees that autocorrelation of the EE sampler decays faster than that of PT. We also plot the autocorrelation for PT with highest temperature $T_7 = 10$ in Figure 7(f), which corresponded to the minimal autocorrelation reached by PT among our tuning values, but the local trapping of PT with this parameter setting is severely pronounced.

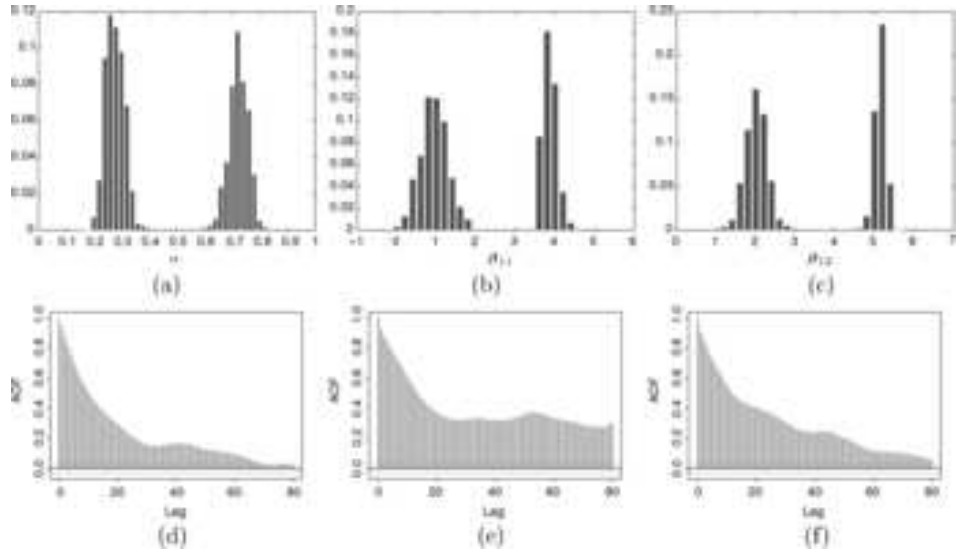

Fig. 7.  *Statistical inference for the mixture exponential regression model. Marginal posterior distribution for* (a) $\alpha$; (b) $\beta_{11}$; *and* (c) $\beta_{12}$ *(the marginal distributions for* $\beta_2$ *are similar to those for* $\beta_1$ *and thus are not plotted here). Autocorrelation plot for the samples from* (d) *EE sampler;* (e) *PT with* $T_7 = 30$; *and* (f) *PT with* $T_7 = 10$.

6.2.  *Motif sampling in biological sequences.*  A central problem in biology is to understand how gene expression is regulated in different cellular processes. One important means for gene regulation is through the interaction between transcription factors (TF's) and their binding sites (viz., the sites that the TF recognizes and binds to in the regulatory regions of the gene). The common pattern of the recognition sites of a TF is called a binding motif, whose identification is the first step for understanding gene regulation. It is very time-consuming to determine these binding sites experimentally, and thus computational methods have been developed in the past two decades for discovering novel motif patterns and TF binding sites.

Some early methods based on site consensus used a progressive alignment procedure [36] to find motifs. A formal statistical model for the position-specific weight matrix (PWM)-based method was described in [21] and a complete Bayesian method was given in [25]. Based on a missing data formulation, the EM algorithm [1, 21] and the Gibbs sampler [20] were employed for motif discovery. The model has been generalized to discover modules of several clustered motifs simultaneously via data augmentation with dynamic programming [41]. See [14] for a recent review.

From a sampling point of view motif discovery poses a great challenge, because it is essentially a combinatorial problem, which makes the samplers very vulnerable to being trapped in numerous local modes. We apply the EE sampler to this problem under a complete Bayesian formulation to test its capability.

6.2.1. *Bayesian formulation of motif sampling.* The goal of motif discovery is to find the binding sites of a common TF (i.e., positions in the sequences that correspond to a common pattern). Let $\mathbf{S}$ denote a set of $M$ (promoter region) sequences, each sequence being a string of four nucleotides, A, C, G or T. The lengths of the sequences are $L_1, L_2, \ldots,$ and $L_M$. For notational ease, let $\mathbf{A} = \{A_{ij}, i = 1, 2, \ldots, M, j = 1, 2, \ldots, L_i\}$ be the indicator array such that $A_{ij} = 1$ if the $j$th position on the $i$th sequence is the starting point for a motif site and $A_{ij} = 0$ otherwise. In the motif sampling problem, a first-order Markov chain is used to model the background sequences; its parameters $\boldsymbol{\theta}_0$ (i.e., the transition probabilities) are estimated from the entire sequence data prior to the motif search. We thus effectively assume that $\boldsymbol{\theta}_0$ is known a priori. Given $\mathbf{A}$, we further denote the aligned motif sites by $\mathbf{S}(\mathbf{A})$, the nonsite background sequences by $\mathbf{S}(\mathbf{A}^c)$ and the total number of motif sites by $|\mathbf{A}|$. The motif width $w$ is treated as known. The common pattern of the motif is modeled by a product multinomial distribution $\boldsymbol{\Theta} = (\boldsymbol{\theta}_1, \boldsymbol{\theta}_2, \ldots, \boldsymbol{\theta}_w)$, where each $\boldsymbol{\theta}_i$ is a probability vector of length 4 for the preferences of the nucleotides $(A, C, G, T)$ in the motif column $i$. See Figure 8 for an illustration of the motif model.

To write down the joint distribution of the complete data and all the parameters, it is assumed a priori that a randomly selected segment of width $w$ has probability $p_0$ to be a motif site ($p_0$ is called the "site abundance" parameter). The joint distribution function has the form

$$P(\mathbf{S}, \mathbf{A}, \boldsymbol{\Theta}, p_0) = P(\mathbf{S}|\mathbf{A}, \boldsymbol{\Theta})P(\mathbf{A}|p_0)\pi(\boldsymbol{\Theta})\pi(p_0)$$

$$(17) \qquad = \frac{P(\mathbf{S}(\mathbf{A})|\mathbf{A}, \boldsymbol{\Theta})}{P(\mathbf{S}(\mathbf{A})|\mathbf{A}, \boldsymbol{\theta}_0)}P(\mathbf{S}|\boldsymbol{\theta}_0)p_0^{|\mathbf{A}|}(1-p_0)^{L-|\mathbf{A}|}\pi(\boldsymbol{\Theta})\pi(p_0)$$

$$\propto \frac{1}{P(\mathbf{S}(\mathbf{A})|\mathbf{A}, \boldsymbol{\theta}_0)} \prod_{i=1}^{w} \boldsymbol{\theta}_i^{\mathbf{c}_i + \boldsymbol{\beta}_i - 1} p_0^{|\mathbf{A}| + a - 1} (1 - p_0)^{L - |\mathbf{A}| + b - 1},$$

where we put a product Dirichlet prior $\pi(\boldsymbol{\Theta})$ with parameter $(\boldsymbol{\beta}_1, \ldots, \boldsymbol{\beta}_w)$ on $\boldsymbol{\Theta}$ (each $\boldsymbol{\beta}_i$ is a length-4 vector corresponding to $A, C, G, T$), and a Beta$(a, b)$ prior $\pi(p_0)$ on $p_0$. In (17), $\mathbf{c}_i$ $(i = 1, 2, \ldots, w)$ is the count vector for the $i$th position of the motif sites $\mathbf{S}(\mathbf{A})$ [e.g., $\mathbf{c}_1 = (c_{1A}, c_{1C}, c_{1G}, c_{1T})$ counts the total number of $A, C, G, T$ in the first position of the motif in all the sequences], the notation $\boldsymbol{\theta}_i^{(\mathbf{c}_i + \boldsymbol{\beta}_i)} = \prod_j \theta_{ij}^{(c_{ij} + \beta_{ij})}$ $(j = A, C, G, T)$, and $L = \sum_{i=1}^{M} L_i$. Since the main goal is to find the motif binding sites given the sequence, $P(\mathbf{A}|\mathbf{S})$ is of primary interest. The parameters ($\boldsymbol{\Theta}$ and $p_0$) can thus be integrated out resulting in the following "collapsed" posterior distribution [14, 24]:

$$
\begin{aligned}
(18) \quad P(\mathbf{A}|\mathbf{S}) \propto{} & \frac{1}{P(\mathbf{S}(\mathbf{A})|\mathbf{A}, \boldsymbol{\theta}_0)} \\
& \times \frac{\Gamma(|\mathbf{A}| + a)\Gamma(L - |\mathbf{A}| + b)}{\Gamma(L + a + b)} \times \prod_{i=1}^{w} \frac{\Gamma(\mathbf{c}_i + \boldsymbol{\beta}_i)}{\Gamma(|\mathbf{A}| + |\boldsymbol{\beta}_i|)},
\end{aligned}
$$

where $\Gamma(\mathbf{c}_i + \boldsymbol{\beta}_i) = \prod_j \Gamma(c_{ij} + \beta_{ij})$ $(j = A, C, G, T)$ and $|\boldsymbol{\beta}_i| = \sum_j \beta_{ij}$.

6.2.2. *The equi-energy motif sampler.* Our target is to sample from $P(\mathbf{A}|\mathbf{S})$. The EE sampler starts from a set of modified distributions,

$$f_i(\mathbf{A}) \propto \exp\left(-\frac{h(\mathbf{A}) \vee H_i}{T_i}\right), \qquad i = 0, 1, \ldots, K,$$

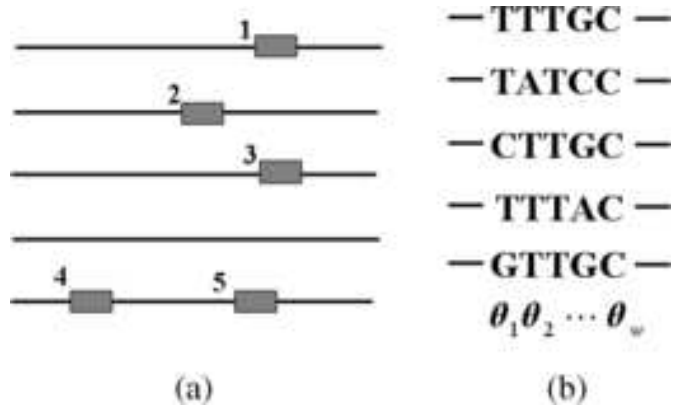

FIG. 8.   (a) *The sequences are viewed as a mixture of motif sites (the rectangles) and background letters. (b) Motif model is represented by PWM, or equivalently, a product multinomial distribution, where we assume that the columns within a motif are independent.*

where $h(\mathbf{A}) = -\log P(\mathbf{A}|\mathbf{S})$ is the energy function for this problem. At the target chain ($i = 0$, $T_0 = 1$), we simply implement the Gibbs sampler, that is, sequentially sample the indicator array $\mathbf{A}$ one position at a time with the rest of the positions fixed. To update other sampling chains at $i = 1, 2, \ldots, K$, given the current sample $\mathbf{A}$ we first estimate the motif pattern $\hat{\boldsymbol{\Theta}}$ by a simple frequency counting with some extra pseudocounts, where the number of pseudocounts increases linearly with the chain index so that at higher-order chains $\hat{\boldsymbol{\Theta}}$ is stretched more toward a uniform weight matrix. Next, we fix $\hat{p}_0 = 1/\bar{L}$, where $\bar{L}$ is the average sequence length. Given $\hat{\boldsymbol{\Theta}}$ and $\hat{p}_0$, we sample each position in the sequences independently to obtain a new indicator array $\mathbf{A}^*$ according to the Bayes rule: Suppose the nucleotides at positions $j$ to $j + w - 1$ in sequence $k$ are $x_1 x_2 \cdots x_w$. We propose $\mathbf{A}^*_{kj} = 1$ with probability

$$q_{kj} = \frac{\hat{p}_0 \prod_{n=1}^{w} \hat{\boldsymbol{\Theta}}_{n x_n}}{\hat{p}_0 \prod_{n=1}^{w} \hat{\boldsymbol{\Theta}}_{n x_n} + (1 - \hat{p}_0) P(x_1 \cdots x_w | \boldsymbol{\theta_0})},$$

where $P(x_1 \cdots x_w | \boldsymbol{\theta_0})$ denotes the probability of generating these nucleotides from the background model. Then $\mathbf{A}^*$ is accepted to be the new indicator array according to the Metropolis–Hastings ratio

$$(19) \qquad r = \frac{f_i(\mathbf{A}^*)}{f_i(\mathbf{A})} \cdot \frac{P(\mathbf{A}|\hat{\boldsymbol{\Theta}}^*)}{P(\mathbf{A}^*|\hat{\boldsymbol{\Theta}})},$$

where $P(\mathbf{A}^*|\hat{\boldsymbol{\Theta}}) = \prod_{k,j} q_{kj}$ denotes the proposal probability of generating the sample $\mathbf{A}^*$ given current $\mathbf{A}$. In addition to the above updating, the EE motif sampler performs the equi-energy jump in each iteration with probability $p_{\mathrm{ee}}$ to help the sampler move freely between different local modes.

6.2.3. *Sampling motifs in "low complexity" sequences.* In the genomes of higher organisms, the presence of long stretches of simple repeats, such as $AAAA\ldots$ or $CGCGCG\ldots$ often makes the motif discovery more difficult, because these repeated patterns are local traps for the algorithms—even the most popular motif finding algorithms based on the Gibbs sampler, such as BioProspector [26] and AlignACE [33], are often trapped to the repeats and miss the true motif pattern. To test whether the EE motif sampler is capable of finding the motifs surrounded by simple repeats, we constructed a set of sequences with the following transition matrix for the background model:

$$(20) \qquad \boldsymbol{\theta}_0 = \begin{bmatrix} 1 - 3\alpha & \alpha & \alpha & \alpha \\ \alpha & 1 - 3\alpha & \alpha & \alpha \\ \alpha & \alpha & 1 - 3\alpha & \alpha \\ \alpha & \alpha & \alpha & 1 - 3\alpha \end{bmatrix},$$

where $\alpha = 0.12$. The data set contained ten sequences, each of length 200 base pairs (i.e., $L_1 = L_2 = \cdots = L_{10} = 200$). Then we independently generated 20 motif sites from the weight matrix whose logo plot [34] is shown in Figure 9(a) and inserted them into the sequences randomly.

The EE motif sampler was applied to this data set with $w = 12$. To set up the energy and temperature ladder, we randomly picked 15 segments of length $w$ in the sequences, treated them as motif sites, and used the corresponding energy value as the upper bound of the energy $H_K$. For the lower bound $H_0$, a rough value was estimated by calculating the energy function (18) for typical motif strength with a reasonable number of true sites. Note that since the EE sampler can adaptively adjust the energy and temperature ladder (see Section 3.3), the bound $H_0$ does not need to be very precise. After some trials $H_0$ was set to be $-50$. We utilized $K + 1 = 5$ chains. The energy ladder $H_j$ was set by a geometric progression between $H_0$ and $H_K$. The temperatures were set by $(H_{j+1} - H_j)/T_j = 5$. The equi-energy jump probability $p_{ee}$ was taken to be 0.1. We ran each chain for 1000 iterations. Our algorithm predicted 18.4 true sites (out of 20) with 1.0 false site on average over ten independent simulations—the EE sampler successfully found the true motif pattern to a large extent. The sequence logo plot of the predicted sites is shown in Figure 9(b), which is very close to the true pattern in Figure 9(a) that generates the motif sites.

The performance of the EE motif sampler was compared with that of BioProspector and AlignACE, where we set the maximum number of motifs to be detected to 3 and each algorithm was repeated ten times for the data set. The motif width was set to be $w = 12$. Both algorithms, however, tend to be trapped in some local modes and output repeats like those shown in Figures 9(c) and (d). One sees from this simple example that the EE sampler, capable of escaping from the numerous local modes, could increase the global searching ability of the sampling algorithm.

6.2.4. *Sampling mixture motif patterns.* If we suspect there are multiple distinct motif patterns in the same set of sequences, one strategy is to introduce more motif matrices, one for each motif type [25]. Alternatively, if we view the different motif patterns as distinct local modes in the sample space, our task is then to design a motif sampler that frequently switches between different modes. This task is almost impossible for the Gibbs sampler, since it can easily get stuck in one local mode (one motif pattern) and have no chance to jump to other patterns in practice. We thus test if the EE sampler can achieve this goal.

In our simulation, we generated 20 sites from each of two different motif models with logo plots in Figures 10(a) and (b) and inserted them randomly into 20 generated background sequences, each of length 100. We thus have 40 motif sites in 20 sequences. The energy ladder in the EE sampler was set by a geometric progression in the range [0, 110] (this is obtained similarly to the previous example). The equi-energy jump probability $p_{ee}$ was set to 0.1. The EE motif sampler used 10 chains; each had 1000 iterations. We recorded the frequency of each position in the sequences being the start of a motif site. Figure 10(c) shows the frequency for the starting positions of the 40 true motif sites (i.e., the probability that each individual motif site was visited by the EE motif sampler). It can be seen that the EE sampler visited both motif patterns with a ratio of 0.4:0.6; by contrast, the Gibbs motif sampler can only visit one pattern in a single run. We further sorted all the positions according to their frequencies of being a motif-site-start in descending order. For any $q$ between 0 and 1, one could accept any position with frequency greater than $q$ as a predicted motif site of the sampler. Thus by decreasing $q$ from 1 to 0, the numbers of both true and false discovered sites increase, which gives the so-called ROC curve (receiver operating characteristic curve) that plots the number of false positive sites versus the number of true positive sites as $q$ varies. Figure 10(d) shows the ROC curves for both the EE motif sampler and the Gibbs motif sampler.

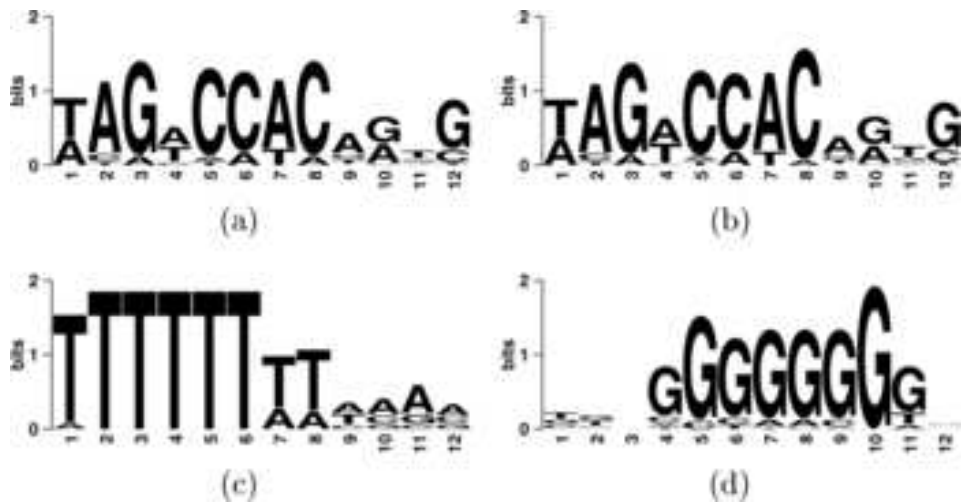

Fig. 9. *Sequence logo plots. (a) The motif pattern that generates the simulated sites. (b) The pattern of the predicted sites by the EE motif sampler. (c) and (d) Repetitive patterns found by BioProspector and AlignACE.*

It can be seen from the figure that for the first 20 true sites, the two samplers showed roughly the same performance. But since the Gibbs sampler missed one mode, the false positive error rate increased dramatically when we further decreased $q$ to include more true sites. By being able to visit both modes (motif patterns), the EE motif sampler, on the other hand, had a very small number of false positive sites until we included 38 true sites, which illustrates that the EE motif sampler successfully distinguished both types of motif patterns from the background sequences.

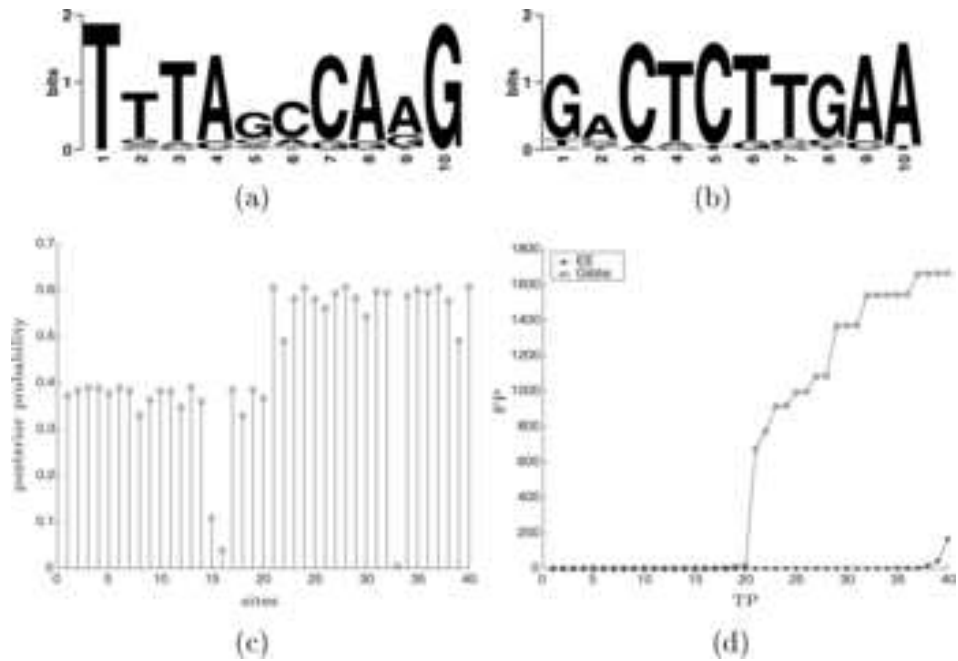

FIG. 10.    (a) and (b) Two different motifs inserted in the sequences. (c) The marginal posterior distribution of the 40 true sites, which reports the frequences that these true site-starting positions are correctly sampled during the iterations. Sites 1–20 are of one motif type; sites 21–40 are of the other motif type. (d) ROC curve comparison between the EE sampler and the Gibbs sampler.

6.3. *The HP model for protein folding.*    Proteins are heteropolymers of 20 types of amino acids. For example, the primary sequence of a protein, cys-ile-leu-lys-glu-met-ser-ile-arg-lys, tells us that this protein is a chain of 10 amino acids linked by strong chemical bonds (peptide bonds) in the backbone. Each different amino acid has a distinct side group that confers different chemical properties to it. Under a normal cellular environment, the side groups form weak chemical bonds (mainly hydrogen bonds) with each other and with the backbone so that the protein can fold into a complex, three-dimensional conformation. Classic experiments such as the ribonuclease refolding experiments [35] suggest that the three-dimensional conformations of most proteins are determined by their primary sequences. The problem of computing the 3D conformation from a primary sequence, known as the protein folding problem, has been a major challenge for biophysicists for over 30 years. This problem is difficult for two reasons. First, there are uncertainties on how to capture accurately all the important physical interactions such as covalent bonds, electrostatic interactions, and interaction between the protein and its surrounding water molecules, and so on, into a single energy function that can be used in minimization and simulation computations. Furthermore, even if the energy function is not in question, we still do not have efficient algorithms for computing the protein conformation from the energy function. For these reasons, biophysicists have developed simplified protein folding models with greatly reduced complexity in both the conformational space and the energy function. It is hoped that the reduced complexity in these simplified models will allow us to deduce insights about the statistical mechanics of the protein folding process through extensive numerical computations.

The HP model is a simplified model that is of great current interest. In this model, there are only two types of amino acids, namely a hydrophilic type (P-type) that is capable of favorable interaction with water molecules, and a hydrophobic type (H-type) that does not interact well with water. Thus, the primary sequence of the length-10 protein in the beginning of this section is simplified to the sequence H-H-H-P-P-H-P-H-P-P. The conformation of the protein chain is specified once the spatial position of each of its amino acids is known. In the HP model space is modeled as a regular lattice in two or three dimensions. Thus the conformation of a length-$k$ protein is a vector $\mathbf{x} = (x_1, x_2, \ldots, x_k)$ where each $x_i$ is a lattice point. Of course, neighbors in the backbone must also be neighbors on the lattice. Figure 11 gives one possible conformation of the above length-10 protein in a two-dimensional lattice model. Just like oil molecules in water, the hydrophobic nature of the H-type amino acids will drive clusters together so as to minimize their exposure to water. Thus we give a favorable energy [i.e., an energy of $\mathcal{E}(x_i, x_j) = -1$] for each pair of H-type amino acids that are not neighbors in the backbone but are lattice neighbors of each other in the conformation [see Figure 11(b)]. All other neighbor pairs are neutral [i.e., having an energy contribution $\mathcal{E}(x_i, x_j) = 0$; see Figure 11(a)]. The energy of the conformation is given by

$$h(\mathbf{x}) = \sum_{|i-j|>1} \mathcal{E}(x_i, x_j).$$

Since its introduction by Lau and Dill [19], the HP model has been extensively studied and found to provide valuable insights [4]. We have applied equi-energy sampling to study the HP model in two dimensions. Here we present some summary results to illustrate how our method can provide information that is not accessible to standard MC methods, such as the relative contribution of the entropy to the conformational distribution, and the phase transition from disorder to order. The detailed implementation of the algorithm, the full description and the physical interpretations of the results are available in a separate paper [16].

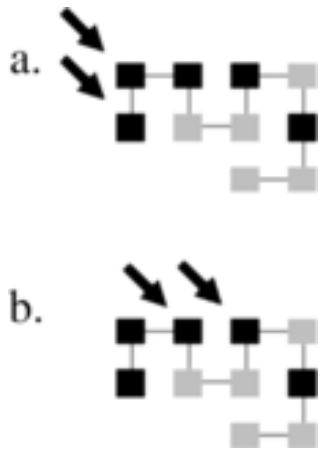

FIG. 11. *One possible conformation of the length-10 proteins in a 2D lattice. H-type and P-type amino acids are represented by black and gray squares, respectively.*

Table 5 presents estimated (normalized) density of states at various energy levels for a protein of length 20: H-P-H-P-P-H-H-P-H-P-P-H-P-H-H-P-P-H-P-H. For this protein there are 83,770,155 possible conformations having energies ranging from 0 (completed unfolded chain) to −9 (folded conformations having the maximum number of nine H-to-H interactions). The estimated and exact relative frequencies of the energies in this collection of conformations are given in the table. The estimates are based on five independent runs each consisting of one million steps where the probability of proposing an equi-energy move is set to be $p_{ee} = 0.1$. It is clear that the method yielded very accurate estimates of the density of states at each energy level even though we have sampled only a small proportion of the population of conformations. The equi-energy move is important for the success of the method: if we eliminate these moves, then the algorithm performs very poorly in estimating the entropy at the low-energy end; for example, the estimated density of state at $E = −9$ becomes $(1.545 \pm 4.539) \times 10^{-12}$, which is four orders of magnitude away from the exact value.

We also use the equi-energy sampler to study phase transition from a disordered state, where the conformational distribution is dominated by the entropy term and the protein is likely to be in an unfolded state with high energy, to an ordered state, where the conformation is likely to be compactly folded structures with low energy. We use the "minimum box size" (BOXSIZE) as a parameter to measure the extent the protein has folded. BOXSIZE is defined as the size of the smallest possible rectangular region containing all the amino acid positions in the conformation. For the 20-length protein, Figure 12 gives a plot of estimated Boltzmann averages of BOXSIZE at ten different temperatures, which are available from a single

Table 5

*The normalized density of states estimated from the EE sampler compared with the actual value*

| Energy | Estimated density of states | Actual value | $t$-value |
|---|---|---|---|
| −9 | $(4.751 \pm 2.087) \times 10^{-8}$ | $4.774 \times 10^{-8}$ | −0.011 |
| −8 | $(1.155 \pm 0.203) \times 10^{-6}$ | $1.146 \times 10^{-6}$ | 0.043 |
| −7 | $(1.452 \pm 0.185) \times 10^{-5}$ | $1.425 \times 10^{-5}$ | 0.144 |
| −6 | $(1.304 \pm 0.189) \times 10^{-4}$ | $1.237 \times 10^{-4}$ | 0.354 |
| −5 | $(9.602 \pm 1.332) \times 10^{-4}$ | $9.200 \times 10^{-4}$ | 0.302 |
| −4 | $(6.365 \pm 0.627) \times 10^{-3}$ | $6.183 \times 10^{-3}$ | 0.291 |
| −3 | $(3.600 \pm 0.228) \times 10^{-2}$ | $3.514 \times 10^{-2}$ | 0.377 |
| −2 | $(1.512 \pm 0.054) \times 10^{-1}$ | $1.489 \times 10^{-1}$ | 0.423 |
| −1 | $(3.758 \pm 0.044) \times 10^{-1}$ | $3.779 \times 10^{-1}$ | −0.474 |
| 0 | $(4.296 \pm 0.071) \times 10^{-1}$ | $4.309 \times 10^{-1}$ | −0.181 |

The $t$-value is defined as the difference between the estimate and the actual value divided by the standard deviation.

run of the equi-energy sampler using five energy ranges. We see that there is a rather sharp transition from order (folded state) to disorder at the temperature range $T = 0.25$ to $T = 1$, with an inversion point around $T = 0.5$. In our energy scale, room temperature corresponds to $T = 0.4$ [16]. Thus at room temperature this length-20 protein will not always assume the minimum energy conformation; rather, it still has a significant probability of being in high-energy, unfolded states. It is clear from this example that the equi-energy sampler is capable of providing estimates of many parameters that are important for the understanding of the protein folding problem from a statistical mechanics and thermodynamics perspective.

**7. Discussion.** We have presented a new Monte Carlo method that is capable of sampling from multiple energy ranges. By using energy truncation and matching temperature and step sizes to the energy range, the algorithm can explore the energy landscape with great flexibility. Most importantly, the algorithm relies on a new type of move—the equi-energy jumps—to reach regions of the sample space that have energy similar to the current state but may be separated by steep energy barriers. The equi-energy jumps provide direct access to these regions that are not reachable with classic Metropolis moves.

We have also explained how to obtain estimates of the density of states and the microcanonical averages from the samples generated by the equi-energy sampler. Because of the duality between temperature-dependent parameters and energy-dependent parameters, the equi-energy sampler thus also provides estimates of expectations under any fixed temperature.

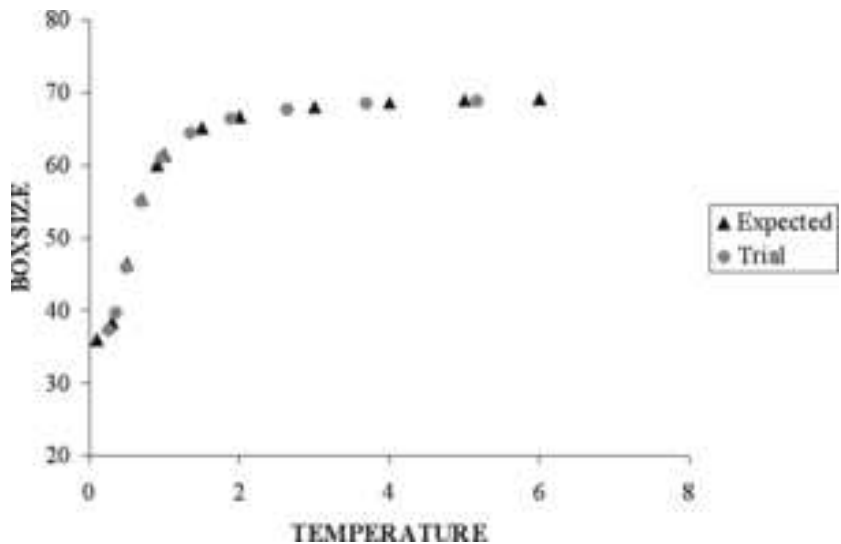

FIG. 12.   *The Boltzmann average of BOXSIZE at different temperatures.*

Our method has connections with two of the most powerful Monte Carlo methods currently in use, namely parallel tempering and multicanonical sampling. In our numerical examples, equi-energy sampling is seen to be more efficient than parallel tempering and it provides estimates of quantities (e.g., entropy) not accessible by parallel tempering. In this paper we have not provided direct numerical comparison of our method with multicanonical sampling. One potential advantage over multicanonical sampling is that the equi-energy jumps allow us to make use of conformations generated in the higher energy ranges to make movement across energy barriers efficiently. Thus equi-energy sampling "remembers" energy-favorable configurations and makes use of them in later sampling steps. By contrast, multicanonical sampling makes use of previous configurations only in terms of the estimated density of states function. In this sense our method has better memory than multicanonical sampling. It is our belief that equi-energy sampling holds promise for improved sampling efficiency as an all-purpose Monte Carlo method, but definitive comparison with both parallel tempering and multicanonical sampling must await future studies.

In addition to its obvious use in statistical physics, our method will also be useful to statistical inference. It offers an effective means to compute posterior expectations and marginal distributions. Furthermore, the method provides direct estimates of conditional expectations given fixed energy levels (the microcanonical averages). Important information on the nature of the likelihood or the posterior distribution, such as multimodality, can be extracted by careful analysis of these conditional averages. For instance, in Section 4 we see that the equi-energy sampler running on (12) reveals that there is a change point in the density of states as well as the microcanonical averages [see Figures 6(c) and (d)], which clearly indicates the multimodality of the underlying distribution. We believe that the design of methods for inferring the properties of the likelihood surface or the posterior density, based on the output of the equi-energy sampler, will be a fruitful topic for future investigations.

**Acknowledgment.** The authors thank Jason Oh for computational assistance and helpful discussion. The authors are grateful to the Editor, the Associate Editor and two referees for thoughtful and constructive comments.

S. C. KOU
Q. ZHOU
DEPARTMENT OF STATISTICS
SCIENCE CENTER
HARVARD UNIVERSITY
CAMBRIDGE, MASSACHUSETTS 02138
USA
E-MAIL: kou@stat.harvard.edu
         zhou@stat.harvard.edu

W. H. WONG
DEPARTMENT OF STATISTICS
SEQUOIA HALL
STANFORD UNIVERSITY
STANFORD, CALIFORNIA 94305-4065
USA
E-MAIL: whwong@stanford.edu